\documentclass[12pt,a4paper]{amsart}
\usepackage{latexsym}
\usepackage{amsmath}
\usepackage{amssymb}
\usepackage{epsfig}
\usepackage{amsfonts}
\usepackage{graphicx}
\usepackage{textcomp}

\begin{document}%
\makeatletter%
\newtheorem{definition}{Definition}[section]%
\makeatother%
\newtheorem{theorem}[definition]{Theorem}
\newtheorem{lemma}[definition]{Lemma}
\newtheorem{proposition}[definition]{Proposition}
\newtheorem{examples}[definition]{Examples}
\newtheorem{corollary}[definition]{Corollary}
\def\square{\Box}
\newtheorem{remark}[definition]{Remark}
\newtheorem{remarks}[definition]{Remarks}
\newtheorem{exercise}[definition]{Exercise}
\newtheorem{example}[definition]{Example}
\newtheorem{observation}[definition]{Observation}
\newtheorem{observations}[definition]{Observations}
\newtheorem{algorithm}[definition]{Algorithm}
\newtheorem{criterion}[definition]{Criterion}
\newtheorem{algcrit}[definition]{Algorithm and criterion}%
\newenvironment{prf}[1]{\trivlist\item[\hskip \labelsep{\it#1.\hspace*{.3em}}]}{~\hspace{\fill}~$\square$\endtrivlist}

\title{Ruled quartic surfaces, models and classification}
\author{Irene Polo-Blanco, Marius van der Put and Jaap Top}
\date{\today}
\address{Department Matesco, University of Cantabria, Avda. Castros s/n, 39005
Santander, Spain.}
\address{IWI-RuG, University of Groningen, Nijenborgh 9, 9747 AG Groningen, 
the Netherlands.}
\email{irene.polo@unican.es, mvdput@math.rug.nl, j.top@rug.nl}

\begin{abstract} New historical aspects of the classification, by Cayley and Cremona,
of ruled quartic surfaces and the relation to string models and plaster models  
are presented. In a `modern' treatment of the classification of ruled quartic surfaces the
classical one is corrected and completed. A conceptual proof is presented of
a result of Rohn concerning curves in $\mathbb{P}^1\times \mathbb{P}^1$ of 
bi-degree $(2,2)$. The string models of Series XIII (of some ruled quartic surfaces) 
are based on Rohn's result.     
\end{abstract}

\maketitle

\section*{Motivation and History}The collection of string models of ruled quartic surfaces, 
present at some mathematical institutes (for instance at
the department of mathematics in Groningen) is the direct motivation for this paper. 
This Series XIII, produced by Martin Schilling in 1886,  is based upon a paper of  
K.~Rohn \cite{Rohn2} containing a classification of ruled quartic surfaces over $\mathbb{C}$ 
and $\mathbb{R}$.  Some authors before Rohn  (e.g., M. Chasles \cite{Chas}, 
A. Cayley \cite{Cay}, L. Cremona \cite{Cre}, R. Sturm \cite{Sturm}, G. Salmon \cite{Sal}) and many
after his time (e.g., B.C. Wong \cite{Wong}, H. Mohrmann \cite{Mohr}, W.Fr. Meyer \cite{Mey}, 
W.L. Edge \cite{Ed}, O. Bottema \cite{Bot}, T. Urabe \cite{Ura2}) have contributed to this 
beautiful topic of 19th century geometry.\\ 

Cremona classified the ruled quartic surfaces in 12 types. He states in \cite{Cre} that 
Cayley produced  8 of these without revealing his method. However, Cayley's third memoir on 
this subject \cite{Cay} was written earlier the same year 1868 and contains 10 types. 
In an addition to this memoir (May 18, 1869), Cayley gives the comparison between his 
own classification and the one by Cremona and makes it clear what the two types he missed are.  
The method of Cayley consists of taking three curves in $\mathbb{P}^3$ and to consider 
the ruled surface $S$ which is the union of the lines meeting all three curves. 
Using  a formula for the degree of $S$, he now computes possibilities of ruled quartic surfaces. 
The expression `the six coordinates of a line' in Cayley's work indicates that the Grassmann 
variety $Gr(2,4)$ of the lines in $\mathbb{P}^3$ plays a role. The work of Cayley contains 
also explicit calculations for reciprocal surfaces (see below). \\ 

The results of Cremona can be explained as follows. Let $S\subset \mathbb{P}^3$ be a ruled 
quartic surface (reduced, irreducible and defined over $\mathbb{C}$).  The fact that through a 
general point of $S$ there is only one line of $S$ is tacitly assumed (compare 
Lemma~\ref{1.2.2}). 
The locus $D$ of the points on $S$ through which there are at least two lines of $S$ 
(in the 1-parameter family) is called the `double curve'. Cremona states that $D$ is indeed a 
curve (hereby excluding cones)  and has `in general' degree 3. 

(We note that $D$ need not 
coincide  with the singular locus of  $S$,  that $D$ can also have degree 2 
(see {\it Number}~$15$ in Subsection~\ref{s1.6}) and that $D=\emptyset$ is possible; 
compare, for example, Corollary~\ref{1.2.7}(2).) 

Two intersecting lines on $S$ determine a plane. The collection $\check{D}$ of all these planes is 
called the `bitangent developable'. This 1-dimensional family (assuming  $D\neq \emptyset$) 
can be seen as a curve in the dual projective space.   The genus of $S$ is defined as the 
genus of the (irreducible, singular) curve $H\cap S$ of degree 4, where $H$ is a general plane.  
Cremona states that the genus can only be $0$ or $1$.  Missing is the nontrivial argument 
showing that  genus 2 is  impossible (see Observation~\ref{1.2.obs} and \cite[Proposition 2.6]{Ura2}). 
Cremona classifies $S$ according to the nature  (degrees and multiplicities of the irreducible 
components) of the curves $D$ and $\check{D}$ (and in one case a relation between $D$ and 
$\check{D}$). He obtains his list of possibilities via the following construction:\\ 

Consider a tuple $(C_1,C_2,f)$ consisting of two conics $C_1, C_2\subset \mathbb{P}^3$, 
not in the same plane, and an isomorphism  $f:C_1\rightarrow C_2$. This defines a ruled 
surface $S$ which is the union of the lines through the pairs of points $\{c_1,f(c_1)\}$ 
with $c_1\in C_1$. In the general case, the line $H_1\cap H_2$, where $C_i$ lies in the 
plane $H_i$ for  $i=1,2$, intersects $C_1$ in two points $p_1\neq q_1$ and intersects $C_2$ 
in two points $p_2\neq q_2$. Now $H_2\cap S$ is the union of the conic $C_2$ and the two lines 
through the pairs of points $(p_1,f(p_1))$ and $(q_1,f(q_1))$. Thus  $S$ is an irreducible ruled 
surface of  degree 4.  Moreover, the two lines intersect in a point of the `double curve' and 
$H_2$ is a `bitangent plane', i.e., a point on the `bitangent developable.'  The same holds 
of course for $H_1$. Cremona's examples are obtained by varying and degenerating 
$C_1,C_2,f$.  His assertion to have found all types in this way is not correct since some 
ruled quartic surfaces are only obtained from a line and a curve of degree 3.  
However by   including `reciprocal surfaces'  and maybe stretching the meaning of `degeneration' 
some of the latter surfaces can be obtained.   

\bigskip

The approach of Cayley (and of Rohn) has the classical name ``analytic geometry'', indicating
the use of coordinates and algebraic operations with formulas. In contrast,
Cremona's (and Sturm's) approach is purely ``synthetic''. As a consequence,
Cremona's paper is difficult to read and it is hard to verify the results.

\bigskip

Bottema \cite[p. 349]{Bot} remarks that Rohn claims to have discovered a type overlooked by 
his predecessors. Indeed, on p.~147 of Rohn's paper \cite{Rohn1}, there is an explicit equation 
and the remark in a footnote: ``this ruled surface is not mentioned by Cremona in his treatise''. 
However, it is easy to verify that Rohn's equation (in the homogeneous coordinates 
$x,y,z,w$)
\[wx^2(x+3Nz)+F_4(x,y)=0,  F_4\mbox{ a binary quartic},  N\mbox{ a constant,} \]
does \textit{not} define a ruled surface, since a general point $(a:b:c:1)$ on it is not contained 
in any line of the surface. Actually, Rohn's geometric construction is valid but his formula 
happens to be mistaken. The construction gives indeed a case which is not  explicitly 
mentioned by Cremona. However it can be interpreted as hiding in Cremona's species 10.  
Pascal's well written Repertorium reviews the classification of Cayley and Cremona, 
\cite{Pas}, XII, \S 10. Here Rohn's extra case reappears on p.~338-339 with the same correct 
geometric construction and another mistake in the formula. In the classification of the present 
paper Rohn's example is {\it Number} 5.

\bigskip

There are also critical comments by R. Sturm to the list of Cremona. Moreover, some of the 
12 species of Cremona contain surfaces of a rather different nature, as we will see in 
Section~\ref{section3}.

The classification of ruled quartic surfaces in the book of W.L. Edge \cite{Ed} is identical 
with the one of  Cremona. Two methods are developed there. The first one classifies the curves 
(irreducible and of degree 4), corresponding to ruled quartic surfaces,  in the Grassmann variety 
$Gr(2,4)$ (parametrizing the 2-dimensional subspaces of a  4-dimensional vector, or, equivalently, 
the lines in $\mathbb{P}^3$). The second method obtains the ruled quartic surfaces in 
$\mathbb{P}^3$ as projections of certain ruled quartic surfaces in $\mathbb{P}^5$ or 
$\mathbb{P}^4$. This is related to a paper of C. Segre \cite{Seg} and to a paper by 
Swinnerton--Dyer\cite{Swin}. \\  

In the thesis of Wong \cite{Wong}, a rational morphism $\mathbb{P}^3 \cdots\to  Gr(2,4)$, 
associated to the  classical `tetrahedral complex', is considered.  Certain plane curves in 
$\mathbb{P}^3$ of degree 2 and 3  have as images in $Gr(2,4)$ curves of degree 4 and 
correspond therefore to ruled quartic surfaces. It is claimed in this thesis  that every `species' 
in Cremona's list can be obtained in this way.\\ 

For other details on the  early history of the subject we refer to the contribution of  
W.Fr.~Meyer in \cite{Mey}. \\

A ruled surface  in modern terminology (see \cite[Section V.2.]{Har}), is a morphism 
$Z\rightarrow C$ of a smooth projective surface $Z$ to a smooth curve $C$ such that all 
fibres are isomorphic to $\mathbb{P}^1$. A classical ruled surface $S\subset \mathbb{P}^3$ is 
obtained as the image of a suitable morphism $Z\rightarrow\mathbb{P}^3$. This method and
 the papers of T. Urabe \cite{Ura1}, \cite{Ura2}, \cite{Ura3}, \cite{Ura4} may lead to a modern
 classification of ruled quartic surfaces (including moduli). We note, in passing, that  Urabe's 
important work concerns the discovery of new aspects in the classification of
the singularities of quartic  surfaces in $\mathbb{P}^3$ (which, generally, are not ruled) and
 their relation to Dynkin diagrams. \\

One aim of the present paper is to give a modern treatment of Rohn's paper \cite{Rohn2},
 namely the`symmetrization' of curves in $\mathbb{P}^1\times \mathbb{P}^1$ of bi-degree  
$(2,2)$ and the classification of some ruled quartic surfaces over $\mathbb{R}$. 
The latter is used to obtain explicit equations explaining the visual features of the models of
 Series XIII.\\

The possibilities for the 1-dimensional part of the singular locus of a ruled surface $S$  can
 be read off from the intersection of $S$ with a general plane. This leads to the elegant
 elementary treatment of ruled cubic surfaces in Dolgachev's book \cite{Dol}.
 As a didactical step towards ruled quartic surfaces, we present here another method valid
 over any base field and obtain the three types of ruled cubic surfaces over $\mathbb{R}$.

\bigskip

The other aim of this paper is to present a classification of the quartic ruled surfaces,
such that each class is determined by discrete data and the surfaces belonging to
a given class give rise to a connected moduli space. This leads to  $29$ cases. 
A combination of the following methods leads to this classification.\\
(1). Deriving some properties of the curves $C$ of degree 4 
(corresponding to ruled quartic surfaces) lying on the Grassmann variety  
$Gr(2,4)\subset \mathbb{P}^5$.\\
(2). Determining the possibilities for the singular locus of a quartic ruled surface.\\
(3). The normalization $C^{norm}$ of $C$ carries a vector bundle of rank two. 
In case the genus of $C^{norm}$ is 0, there are two possibilities for this vector bundle.
Two `generating' meromorphic sections of this vector bundle are brought in some
standard form, by some linear base changes. This  leads to explicit equations for the 
corresponding  quartic ruled surfaces.\\
(4). Classification of the position of $C$ w.r.t. the tangent spaces of $Gr(2,4)$.\\

Part (4) is in fact one of the two methods of \cite{Ed} in deriving Cremona's list. Although
we could not verify this in detail because of a certain vagueness in Edge's arguments,
the results agree with our computations.\\

Before giving Cremona's list we need to explain to notion of `{\it reciprocal surface}' or 
`{\it dual surface}' in more modern terms, of a surface $S\subset \mathbb{P}^3$. It  is obtained
 by considering all tangent planes at the nonsingular points of $S$. Each tangent
plane is a point in the dual projective space $\check{\mathbb{P}^3}$. The Zariski closure
 of all these points is the dual surface $\check{S}\subset \check{\mathbb{P}^3}$.
  In the case that $S$ is ruled, also $\check{S}$ is ruled and has the same degree as $S$.
 Moreover, the `double line' of $\check{S}$ can be seen to be the `bitangent developable'.
 Cremona's list shows that the species 3 and 4 are dual, as well as the species 7 and 8. 
The other species are `selfdual'. In the table one has to give `double curve' the interpretation 
`singular locus'. The genus $g$ of the surface is 0, except for Cremona 11, 12 where it is 1.

We adopt a {\it notation of Cayley}, namely the expression $d^m$ stands for an irreducible  
component of the singular locus of degree $d$ and with multiplicity $m$. 
The difference between Cremona 6 and Cremona 11 is that the two lines intersect in the first case
 and are skew in the second one. The difference between Cremona 9 and Cremona 10 is
  somewhat subtle. In case 10 the bitangent planes are the planes containing the singular line,
 denoted by $1$. In case 9, the bitangent planes are the planes containing another line,
 denoted by $1'$. 

\begin{center}        
$  \begin{array}{| c | c | l|| c |c | l|}       \hline         
\mbox{Double} & \mbox{Double curve} & \mbox{Cremona} &         
\mbox{Double} & \mbox{Double curve}  & \mbox{Cremona}\\ 
       \mbox{curve} &  \mbox{recip. surface} & \mbox{type} & \mbox{curve} & 
 \mbox{recip. surface} & \mbox{type} \\        
\hline       3^2&3^2&1&3^2&1^3 &7\\       \hline  
      2^2,1^2&2^2,1^2&2&1^3&3^2 &8\\       \hline  
      1^3&2^2,1^2&3&1^3&{1'}^3&9\\       \hline       
  2^2,1^2& 1^3&4& 1^3& 1^3&10\\       \hline     
   1^2,1^2,1^2&1^2,1^2,1^2&5&1^2,1^2&1^2,1^2&11,\ g=1\\       \hline    
   1^2,1^2\mbox{ int} &1^2,1^2\mbox{ int}&6 &1^2&1^2&12,\ g=1\\       \hline    
   \end{array} $   \end{center}
   
   \smallskip

\begin{center}
\scalebox{.42}{\includegraphics{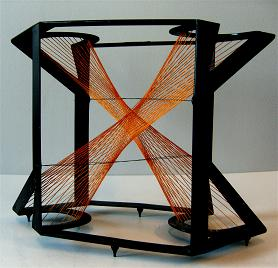}}
\scalebox{.45}{\includegraphics{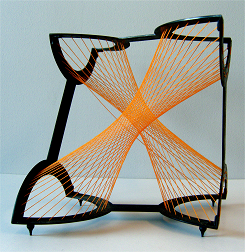}}
\scalebox{.41}{\includegraphics{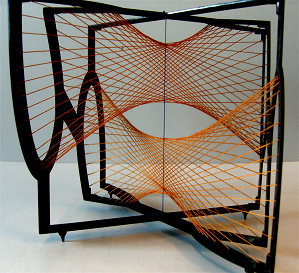}}\\
Some models from Series~XIII: nr.\  1, 4, and 5.

\end{center}

\section{Curves on the Grassmann variety $Gr(2,4)$}

\subsection{Properties of the Grassmann variety}
Let $V$ be a vector space of dimension 4 over the (algebraically closed) field $K$.
 The lines in theprojective space $\mathbb{P}(V)\cong \mathbb{P}^3$ are points of the 
Grassmann variety $Gr:=Gr(2,V)=Gr(2,4)$ and the natural way to study a ruled surface 
$S\subset \mathbb{P}(V)$ is to consider the set of the lines on $S$ as subset of $Gr$.
  We briefly define $Gr$ and summarize its main properties.
  
For notational convenience we fix a basis $e_1,e_2,e_3,e_4$ of $V$ and we identify the
 exterior power$\Lambda ^4 V$ with $K$ by  $e_1\wedge \cdots \wedge e_4 \mapsto 1$. 
The obvious symmetric bilinear map 
$\Lambda ^2V\times \Lambda ^2V\rightarrow \Lambda ^4V=K$ is nondegenerate. 
A line in $\mathbb{P}(V)$ correspond to a plane $W\subset V$,  
a line $\Lambda ^2W\subset  \Lambda ^2V$ and to a point in 
$\mathbb{P}(\Lambda ^2V)\cong \mathbb{P}^5$. If $W$ has basis $v_1,v_2$, then 
$w=v_1\wedge v_2$ is a basis vector for $\Lambda ^2W$ and $\overline{w}:=Kw$ is this point 
of $\mathbb{P}(\Lambda ^2V)$. By definition $Gr=Gr(2,V)\subset \mathbb{P}(\Lambda ^2V)$ consists of all these 
points. Now $\overline{w}$ (with $w\in \Lambda ^2V,\ w\neq 0$) belongs to $Gr$ if and only if $w$ is 
decomposable, i.e., has the form $v_1\wedge v_2$.  The latter is equivalent to $w\wedge w=0$. 
We use the six elements $e_{ij}:=e_i\wedge e_j,\ i<j$ as basis for $\Lambda ^2V$ and write an 
element of this vector spaceas $\sum _{i<j} p_{ij}e_{ij}$. The $p_{ij}$ are called the 
Pl\"ucker coordinates. They also serve as homogeneous coordinates for 
$\mathbb{P}(\Lambda ^2V)$. One finds that $Gr$ is the nondegenerate quadric given by the 
equation $p_{12}p_{34}-p_{13}p_{24}+p_{14}p_{23}=0$. For notational purposes and for 
convenience of the reader we recall the following.\\

{\bf List of properties of $Gr$} (of importance for our purposes).  
\begin{enumerate}   
\item[(i)] $p_0,\ \ell _0, h_0$ are a point, a line and a plane of $\mathbb{P}(V)=\mathbb{P}^3$.  
One identifies $p_0$ with a $\overline{v}_0,\ v_0\in V,\ v_0\neq 0$ and $\ell _0$ with a  
$\overline{w}_0,\ w_0\in \Lambda ^2V,\ w_0\wedge w_0=0$. 
{\it We note that $\overline{w}\in \mathbb{P}(\Lambda ^2V)$ with $w\wedge w=0$ is both seen as a 
point of $Gr$ and as a line in $\mathbb{P}(V)$}.      
\item[(ii)] Two lines $\overline{w}_1,\overline{w}_2$ of $\mathbb{P}(V)$ intersect if and only if  
$w_1\wedge w_2=0$.     
 \item[(iii)] Every hyperplane of $\mathbb{P}(\Lambda ^2V)$ has the form    
$\{\overline{z}| \ w\wedge z=0\}$ with $w \in \Lambda ^2V, \ w\neq 0$ and unique    
$\overline{w}$.  If $w$ is indecomposable, then the intersection of the hyperplane   
with $Gr$ is a nondegenerate quadric.       
If $w$ is decomposable, i.e., $\overline{w}=\overline{w}_0\in Gr$ and correspond to the line
 $\ell _0$,   then the hyperplane is the  tangent plane $T_{Gr,\overline{w}_0}$ of $Gr$ at
 $\overline{w}_0$.   The intersection $T_{Gr,\overline{w}_0}\cap Gr$ is singular and can be 
identified with the cone   in $\mathbb{P}^4$ over a nonsingular quadric in $\mathbb{P}^3$.       
This intersection identifies with $\sigma _1(\ell _0):=$  the collection of all lines $\ell$ with    
$\ell \cap \ell _0\neq \emptyset$.  Consider for example  $w_0=e_1\wedge e_2$. 
This intersection    is now $\{\{p_{ij}\}|\ p_{34}=0,\ -p_{13}p_{24}+p_{14}p_{23}=0\}$. 
The vertex $\overline{w}_0$ of this     cone is its only singular point.        
\item[(iv)] $\sigma _2(p_0):=$ the collection of all lines through $p_0$; this is a 2-plane in $Gr$.     
 Indeed, take $p_0=\overline{e}_1$. Then 
\[\sigma _2(p_0)=     
\left\{\overline{\sum _{1<j\leq 4}p_{1j}e_1\wedge e_j}| \mbox{ no relations}\right\} \]    
(called a $\omega$-plane in \cite{Ed}).    \item[(v)] $\sigma _{1,1}(h_0):=$ the collection of all lines 
in the plane $h_0$. This is a 2-plane    in $Gr$. Indeed, take $h_0=\overline{<e_1,e_2,e_3>}$. 
Then this collection identifies with    
$\{ \overline{\sum _{1\leq i<j\leq 3}p_{ij}e_i\wedge e_j}|\mbox{ no relations }\}$ 
(called a    $\rho$-plane in \cite{Ed}).        
\item[(vi)] $\sigma _{2,1}(p_0,h_0):=$ the collection of all lines in $h_0$ through $p_0$. This is    
a line on $Gr$. Indeed, take $h_0=\overline{<e_1,e_2,e_3>},\ p_0=\overline{e}_1$. Then this    
collection identifies with     
$\{\overline{\sum _{j=2,3}p_{1j}e_1\wedge e_j}|\mbox{ no relations}\}$.     
 \item[(vii)] Every plane in $Gr$ has the form $\sigma _2 (p_0)$ or $\sigma _{1,1}(h_0)$.   
 Every line in $Gr$ has the form $\sigma _{2,1}(p_0,h_0)$ and is thus the intersection of     
a (uniquely determined) pair of  2-planes in $Gr$ of different type.         
\item[(viii)]  The are three types of projective subspaces $P\subset \mathbb{P}(\Lambda ^2V)$ of 
dimension 3    with respect to their relation with $Gr$,  namely:\\    
(a) $Gr\cap P$ is a nondegenerate quartic surface. The equations of $P$ are $p_{12}=p_{34}=0$ 
for a suitable     basis of $V$. Moreover $P$ lies in precisely two tangent space, 
namely  $T_{Gr,\ \overline{e_{12}}}$ and    $T_{Gr,\ \overline{e_{34}}}$.\\    
(b) $Gr\cap P$ is an irreducible  degenerate quartic surface. The equations of $P$     
are $p_{34}=p_{13}+p_{24}=0$  for a     suitable basis of $V$. Now $P$ lies on only one 
tangent space, namely $T_{Gr,\ \overline{e_{12}}}$ and    $Gr\cap P$ is the cone 
$p_{13}^2+p_{14}p_{23}=p_{34}=p_{13}+p_{24}=0$ over the quadric curve    
$p_{13}^2+p_{14}p_{23}=0$.\\    
(c) $Gr\cap P$ is reducible. The equations for $P$ are $p_{12}=p_{13}=0$ for a suitable 
basis of $V$.    Further, $Gr\cap P$ is the union of the  planes $p_{14}=0$ and $p_{23}=0$.   
\hfill $\square$ 
\end{enumerate}        
\subsection{Ruled surfaces and curves on $Gr$}   
\begin{lemma}\label{1.2.1} {\rm (1)}. Let $C\subset Gr$ be an irreducible curve of degree    
$d\geq 2$, not lying in some 2-plane $\sigma _2 (p_0)$.    
Then $\tilde{S}:=\{(\overline{w},\overline{v})\in C\times \mathbb{P}(V)|\ w\wedge v=0\}$   
is an irreducible variety of dimension 2. Its image $S$ under the projection map   
$pr_2 : \tilde{S}\rightarrow \mathbb{P}(V)$ is an irreducible surface of degree $e$.    
Suppose that through a general point of $S$ there are $f$ lines $\overline{w}\in C$.     
Then $d=e\cdot f$.\\    
{\rm (2)}. Let $P(C)$ denote the smallest projective subspace of    
$\mathbb{P}(\Lambda ^2V)$, containing $C$.  If $d\geq 3$ and $S$ is not a cone,     
a plane or a quadric, then  $\dim P(C)\geq 3$.   
\end{lemma}   

\begin{proof} (1). We note that $C\subset  \sigma _2(p_0)$ is not interesting since  
then $S$ is a cone.  The fibers of $pr_1:\tilde{S}\rightarrow C$ are lines in   
$\mathbb{P}(V)$ and the fibers of $pr_2:\tilde{S}\rightarrow S$ are finite. Thus $S$   
is an irreducible  ruled surface of some degree $e$.  A general line $\overline{w}_0$ in   
$\mathbb{P}(V)$ intersects $S$ in $e$ points. Through each of these $e$ points there    
are $f$ lines $\overline{w}\in C$. Thus the intersection of $C$ with the general    
hyperplane  $\{\overline{w}\in \mathbb{P}(\Lambda ^2V)|\ w\wedge w_0=0\}$ 
consists of $e\cdot f$ points and  therefore $d=e\cdot f$.\\  
(2).  Since $d>1$, one has $\dim P(C)>1$. Suppose that  $\dim P(C) =2$.   
If  $P(C)\subset Gr$, then either $P(C)$ is  a $\sigma _2(p_0)$ and $S$ is a cone, or   
$P(C)$ is a $\sigma _{1,1}(h_0)$ and $S$ is the plane $h_0$. If $P(C)\not \subset Gr$,   
 then $C\subset P(C)\cap Gr$ is a curve of degree at most 2 and $S$ is a plane     
or a quadric. Hence $\dim P(C)\geq 3$. 
\end{proof}     

 In the sequel we consider ruled surfaces (reduced, irreducible) $S\subset \mathbb{P}(V)$ of 
some degree $d\geq 3$ which are not cones. One associates to $S$ the subset  
$\tilde{C}$ of $Gr$ corresponding to the lines on $S$.           

\begin{lemma}\label{1.2.2} $\tilde{C}$ is the union of an irreducible curve $C$     
(not lying in some 2-plane $\sigma _2(p_0)$) of degree $d$ and a finite,  possibly empty, set. 
Moreover, through a general point of $S$ there is one line of the surface.      
\end{lemma}    

 \begin{proof}  Consider the affine open part of $Gr$ given by $p_{12}\neq 0$. The points of this      
affine part, actually $\cong \mathbb{A}^4$, can uniquely be written as planes in $V$ with basis       
$e_1+ae_3+be_4,\ e_2+ce_3+de_4$  and correspond to the  vectors      
\[e_{12}+ce_{13}+de_{14}+ae_{23}+be_{24} +  (ad-bc)e_{34}\ .\]      
Let $F(t_1,\dots ,t_4)=0$ be the homogeneous equation of $S$. The intersection of $\tilde{C}$       
with this affine part consists of the  tuples $(a,b,c,d)$ such that $F(s,t,as+ct,bs+dt)=0$ for all        
$(s,t)\neq (0,0)$. Write this expression as a homogeneous form in $s,t$ and coefficients   
 polynomials in $a,b,c,d$. Then the ideal generated by these polynomials in $a,b,c,d$ defines the   
intersection of $\tilde{C}$ with this affine part of $Gr$. 
Thus $\tilde{C}$ is Zariski closed. 
         
Clearly   $\tilde{C}$ has dimension 1 and can be written as the union of irreducible curves 
$C_i,\ i=1, \dots, r$ and a finite set. The image of  the projection 
$\{ (\overline{w},\overline{v})\in C_1\times \mathbb{P}(V) | w\wedge v=0\}
\rightarrow \mathbb{P}(V)$ is  a ruled surface contained in $S$. Since $S$ is irreducible, 
the image is $S$. If $r\geq 2$, then, for through a point $\overline{v}$ of a line   
$\overline{w}_2\in C_2,\ \overline{w}_2\not \in C_1$ passes a line  $\overline{w}_1\in C_1$. 
Hence $w_1\wedge w_2=0$ for all $\overline{w_1}\in C_1$ and thus $w\wedge w_2=0$ for all 
$\overline{w}\in P(C_1)$. By symmetry  $w_1\wedge w_2=0$ for all 
$\overline{w}_1\in P(C_1),\ \overline{w}_2\in P(C_2)$. 
Since the symmetric bilinear form $(w_1,w_2)\mapsto w_1\wedge w_2$ on $\Lambda ^2V$ 
is not degenerate, one obtains a contradiction by compairing  dimensions: 
$\dim P(C_1)\geq 3,\ \dim P(C_2)\geq 3,\ \dim \mathbb{P}(\Lambda ^2V)=5$. 
We conclude that the $f$ of Lemma~\ref{1.2.1} is 1 and that the degree of $C$ is $d$. 
\end{proof}

\begin{lemma}\label{1.2.3} Let $\overline{w}_0\in \tilde{C}\setminus C$, then $C$ lies in the tangent space of 
$Gr$ at $\overline{w}_0$.  In other words, the line $\overline{w}_0$ intersects every line on $S$, 
belonging to $C$. 
\end{lemma}

\begin{proof}  If the tangent space at $\overline{w}_0$  does not contain $C$, then the intersection
 $C\cap T_{Gr,\ \overline{w}_0}$ consists of $d$ points, counted with multiplicity. 
Thus the line $\overline{w}_0$ on $S$ intersects $d$ lines of $S$, corresponding to points of $C$. 
Let $H\subset \mathbb{P}(V)$ be a plane through $\overline{w}_0$. 
The intersection $H\cap S$ consists of $\overline{w}_0$ and a curve $\Gamma$ of degree $d-1$.  
Therefore $\Gamma \cap \overline{w}_0$ consists of $d-1$ points (counted with multiplicity),
 instead of the $d$ points that we expect. This contradiction proves the lemma.
\end{proof}

\begin{definition}. {\rm  The lines on $S$ corresponding to the points of $\tilde{C}\setminus C$ 
will be called here {\it isolated lines}. A line $\overline{w}_1$ on $S$ is, classically, 
called a {\it directrix} if $\overline{w}_1$ meets every line $\overline{w}_2$ with 
$\overline{w}_2\in C$. Thus an isolated line is a directrix. It is also possible that a line 
$\overline{w}_1\in C$ is a directrix. The classical concept of `{\it double curve}' on $S$ is, 
according to \cite{Ed}, p. 8, (the Zariski closure of) the set of points on $S$ lying on at least two, 
not isolated, lines of $S$. }
\end{definition}   

\begin{observation}\label{1.2.obs} Let $C\subset P=\mathbb{P}^r, r\geq 3$ be an irreducible curve of degree 4 
and such that $C$ does not lie in a proper subspace of $P$. Let $g\leq 2$ be the genus of the 
normalization $n:C^{norm}\rightarrow C$. Then one of the following holds:\\  
(1)  $g=0,\ r=4$, $C$ is the, nonsingular, rational normal quartic.\\ 
(2)  $g=0,\ r=3$, $C$ is nonsingular or has one singular point which is a node or ordinary cusp.\\ 
(3) $g=1,\ r=3$ and  $C$ is nonsingular.\\ 
Moreover, if $C$ lies on a quadratic cone in $\mathbb{P}^3$, then $g=1$ or  $g=0$ and $C$ 
has a singular point. 
\end{observation} 

\begin{proof}  Let $I$ be the sheaf of ideals of $C$. The exact sequence 
$0\rightarrow I\otimes O_P(1)\rightarrow O_P(1)\rightarrow \mathcal{L}\rightarrow 0$ with
$\mathcal{L}=O_P(1)\otimes O_C$ and the minimality of $r$ implies that
$H^0(O_P(1))\rightarrow H^0(\mathcal{L})$ is injective and thus $1+r\leq \dim H^0(\mathcal{L})$.
Define the skyscraper sheaf $\mathcal{Q}$ on $C$ by the exact sequence of sheaves on 
$C$, $0\rightarrow \mathcal{L}\rightarrow n_*n^*\mathcal{L}\rightarrow \mathcal{Q}\rightarrow 0$.
Denoting $\dim H^i$ by $h^i$, one finds  
\[ 4\leq 1+r\leq h^0(C,\mathcal{L})\leq h^0(C^{norm},
n^*\mathcal{L})=1-g+4+\dim H^1(C^{norm}, n^*\mathcal{L}) .\]
Now $H^1(C^{norm},n^*\mathcal{L})=0$, since the degree of $n^*\mathcal{L}$ is 4 and $g\leq 2$.
Thus $g=2$ is not possible. For $g=1$, one has $H^0(C,\mathcal{L})= H^0(C^{norm},n^*\mathcal{L})$ 
and $C=C^{norm}$ since $n^*\mathcal{L}$ is very ample on $C^{norm}$. Let $E$ be an elliptic  
curve with neutral element $e$, then $H^0(E, 4[e])$ has basis $t_1=1,t_2=x,t_3=y,t_4=x^2$ 
(in the standard notation) and $E$ lies on the quadratic cone $t_2^2-t_1t_4=0$. 

For $g=0$, the curves $C\subset P$ are parametrized by polynomials of degree $\leq 4$ in a 
variable $t$. Hence $r\leq 4$. For $r=4$, the only possibility is  
$t\mapsto (1,t,t^2,t^3,t^4)$. For $r=3$, one has the examples:\\ 
$t\mapsto (1,t,t^3,t^4)$ and $C$ is nonsingular,\\
$t\mapsto (1,t^2,t^3,t^4)$ and $C$ has an ordinary cusp,\\
$t\mapsto   (t,t^2,t^3,t^4-1)$ and $C$ has a  node.\\

In general, by intersecting $C$ with planes $H\subset \mathbb{P}^3$, through one singular point 
(or more), one can verify that $C$ has at most one singular point and that such a point can
 only be a node or an ordinary cusp. Finally, if $g=0$ and $C$ is contained in a quadratic cone in 
$\mathbb{P}^3$, then $C$ is singular (see  \cite{Har}, exercise IV, 6.1). According to the examples, 
this singular point can be either a node or a cusp.
We note that the exercises IV, 3.4, 3.6 and II, Example 7.8.6 of \cite{Har} are
closely related to the above reasoning. 
\end{proof}

\begin{corollary}\label{1.2.6} A ruled surface of degree $d\geq 3$ can have at most two isolated lines. 
If $S$ has two isolated lines $\overline{w}_1,\ \overline{w}_2$, then 
$\overline{w}_1\cap \overline{w}_2=\emptyset$.
\end{corollary}

\begin{proof} The first statement follows from $\dim P(C)\geq 3$. If 
$\overline{w}_1\cap \overline{w}_2\neq \emptyset$,
then $C$ lies in $Gr\cap T_{Gr,\ \overline{w}_1}\cap T_{Gr,\overline{w}_2}$. 
According to the list of properties of $Gr$, (viii) part (c), the latter is the union of two planes. 
One of them contains $C$ and this contradicts $\dim P(C)\geq 3$.  
\end{proof}

\begin{corollary}\label{1.2.7} {\rm (1)} A general line of a `general' ruled  surface $S$ of degree $d\geq 3$ 
meets $d-2$ other lines of $S$, corresponding to points of $C$. In particular, 
the `double curve' is not empty. {\rm However:} \\ 
{\rm (2)}  Let $TC\subset \mathbb{P}^3$ be the twisted cubic curve. The equation of the surface 
$S$ consisting of the tangents  of  $TC$  is 
$(t_1t_4-t_2t_3)^2-4(t_1t_3-t_2^2)(t_2t_4-t_3^2)=0$.
The singular locus of $S$ is $TC$ and no two distinct lines of the surface intersect. 
\end{corollary}

\begin{proof} (1) For a general point $\overline{w}_0\in C$, the intersection 
$C\cap T_{Gr,\ \overline{w}_0}$is a positive divisor on $C$ of degree $d$, with support in the 
nonsingular locus of $C$ and $\geq 2[\overline{w}_0]$. For a `general' $S$ the divisor will be 
$2[\overline{w}_0]+\sum _{i=1}^{d-2} [\overline{w}_i]$
with distinct points $\overline{w}_i\in C, \ i=0,\dots , d-2$. Thus $\overline{w}_0$ meets 
precisely $d-2$ other lines corresponding to points of $C$.\\
\noindent (2) Let $t\mapsto (1,t,t^2,t^3)\in \mathbb{P}^3$ be $TC$ in parametrized form. 
The tangent line $\overline{w}_t$ contains the point $(0,1,2t,3t^2)$ and has 
Pl\"ucker coordinates
\[p_{12}=1,\ p_{13}=2t,\ p_{14}=3t^2,\ p_{23}=t^2,\ p_{24}=2t^3,\ p_{34}=t^4 .\]  
This defines the nonsingular curve $C\subset Gr$ corresponding to $S$. 
From $\overline{w}_t\wedge \overline{w}_s=(t-s)^4$ it follows that
the tangent lines do not intersect for $t\neq s$. In other terms 
$T_{Gr,\overline{w}_0}\cap C =4[\overline{w}_0]$ for every 
$\overline{w}_0\in C$.  
\end{proof}

\begin{proposition}\label{1.2.8}  
$\ $\\{\rm (1)} 
$pr_2: \tilde{S}:=\{(\overline{w},\overline{v})\in C\times \mathbb{P}(V)|\ w\wedge v=0\} \rightarrow S$ 
is a birational morphism.  Let $C^{norm}\rightarrow C$ denote the normalization of $C$ 
and let $\tilde{\tilde{S}}=\tilde{S}\times _{C}C^{norm}$ be the pullback  of $\tilde{S}\rightarrow C$.  
Then  $\tilde{\tilde{S}}\rightarrow C^{norm}$ is a ruled surface (in the modern sense) and 
$\tilde{\tilde{S}}\rightarrow S$ is the normalization of $S$.\\
\noindent {\rm (2)} The singular locus of $S$ is purely 1-dimensional or empty.\\
\noindent {\rm (3)} Suppose that the line $\overline{w}$ belongs to the singular locus of $S$ 
and does not correspond to a singular point of $C$. Then $C$ lies in the tangent space of $Gr$ 
at the point $\overline{w}$.
\end{proposition}

\begin{proof} (1) The morphism is finite since $pr^{-1}_2(\overline{v})$ is the finite set of lines 
of $S$ through $\overline{v}\in S$. For a general $\overline{v}\in S$, this set has one element 
and therefore  the degree of $pr_2$ is 1 and so $pr_2$ is birational. 
The fibres of $pr_1:\tilde{S}\rightarrow C$ are isomorphic to $\mathbb{P}^1$ and the same holds 
for the fibres of $\tilde{\tilde{S}}\rightarrow C^{norm}$. Therefore the latter is a ruled surface in the 
modern terminology. Moreover the morphism $\tilde{\tilde{S}}\rightarrow \tilde{S}$ is birational 
and so $\tilde{\tilde{S}}\rightarrow S$ is the normalization. \\
\noindent (2) The local ring of an isolated singular point of $S$ is normal  and will remain a 
singular point of the normalization of  $S$. Since $\tilde{\tilde{S}}$ is smooth, 
$S$ has no isolated singularities. \\  
\noindent (3) The assumption that $\overline{w}$ does not correspond to a singular point of $C$ 
implies that through any point of $\overline{w}$ there are at least two lines corresponding to points 
of $C$ (one of them could be $\overline{w}$ itself). Hence $\overline{w}$ meets every line 
corresponding to a point of $C$ and thus $C\subset T_{Gr,\overline{w}}$.  
\end{proof}

\begin{remarks} {\rm   The  `double curve', as defined above,  is seen, by
Proposition~\ref{1.2.8}, 
to be part of the  singular locus of $S$. The {\it genus of $S$} is defined as the genus of 
$\tilde{\tilde{S}}$ and thus equal to the genus of $C^{norm}$.}\hfill $\square$ 
\end{remarks}    

\begin{lemma}\label{1.2.10} Suppose that $\dim P(C)=3$ and  that $P(C)$ is the intersection of two tangent 
spaces of $Gr$ at points  $\overline{w}_1\neq  \overline{w}_2$. Then the lines 
$\overline{w}_1,\ \overline{w}_2$ do not intersect.  For a suitable choice of the homogeneous 
coordinates $t_1,t_2,t_3,t_4$ of $\mathbb{P}(V)$, the equation $F$ of $S$ is bi--homogeneous 
of degree $(a_1,a_2)$, with $a_1+a_2=d$, in the pairs $t_1,t_2$ and $t_3,t_4$.  
Further $\tilde{C}\setminus C= \{\overline{w}_1,\overline{w}_2\}$. 
  
The lines $\overline{w}_1, \overline{w}_2$ are `directrices'. The singular locus of $S$  
consists of  the lines $\overline{w}_i$ with $a_i>1$ and for each singular point 
$\overline{w}\in C$, the line   $\overline{w}\subset S$.
\end{lemma} 

\begin{proof}  The assumption that the lines $\overline{w}_1,\ \overline{w}_2$ intersect, yields, 
according to (viii) part (c), the contradiction that $C$ lies in a plane. 
Take $w_1=e_{12}$ and $w_2=e_{34}$, then  
$P(C)= T_{Gr,\ \overline{e_{12}}}\cap T_{Gr,\ \overline{e_{34}}}$ is the projective space with 
coordinates $p_{13},p_{14},p_{23},p_{24}$ and $C$ lies on the quadric surface 
$Gr\cap P(C)$ given by $-p_{13}p_{24}+p_{14}p_{23}=0$. Identifying $Gr\cap P(C)$ with 
$\mathbb{P}^1\times \mathbb{P}^1$ leads to $C\subset \mathbb{P}^1\times \mathbb{P}^1$ 
of bi--degree $(a_1,a_2)$ with $a_1+a_2=d$.
    
Consider the rational map $f:\mathbb{P}(V)\cdots\rightarrow \mathbb{P}^1\times \mathbb{P}^1$, 
given by 
\[(t_1,t_2,t_3,t_4)\mapsto ((t_1,t_2),(t_3,t_4)),\] 
which  is defined outside the two lines 
$\overline{w}_1,\overline{w}_2$. The surface $S$ is the Zariski closure of $f^{-1}(C)$ and so 
the equation $F$ of $S$ is bi-homogeneous and coincides with the equation for 
$C\subset \mathbb{P}^1\times \mathbb{P}^1$. The other statements of the lemma are easily 
verified.    
\end{proof}  

\begin{observation}\label{1.2.11} $\dim P(C)=3$ and $P(C)$ in a single tangent space of $Gr$.\\
{\rm  For a suitable basis of $V$ the projective subspace 
$P(C)\subset T_{Gr,\ \overline{e_{12}}}$ is given by the equations $p_{34}=0, p_{13}+p_{24}=0$ 
and  $p_{12},p_{13},p_{14},p_{23}$ are the homogeneous coordinates of $P(C)$. 
Further $Gr\cap P(C)$ is the cone with equation $p_{13}^2+p_{14}p_{23}=0$ with vertex 
$\overline{e_{12}}$. Since  $C$ lies on this cone we have a rational map 
$f: C\cdots \rightarrow E:=\{p_{13}^2+p_{14}p_{23}=0\}$. This map can be identified with the 
rational map $C\cdots\rightarrow \overline{e_{12}}$, given by  
$\overline{w}\mapsto \overline{w}\cap \overline{e_{12}}$. The rational map  
$f$ is a morphism if $\overline{e_{12}}\not \in C$ or if $\overline{e_{12}}\in C$ and this 
is a regular point of $C$. 

In case $\overline{e_{12}}\not \in C$ the morphism $f$ has 
degree $e$. Take two unramified points $e_1,e_2\in E$ and the plane through the corresponding 
two lines through $\overline{e_{12}}$. This plane meets $C$ in $2e$ points. 
Hence $d=2e$. In case $\overline{e_{12}}\in C$ and is not a singular point, the same reasoning 
yields $d-1=2e$. 

It seems difficult to investigate the possibilities for general $d$. 
The cases $d=3$ and $d=4$ will be presented later on. }
\end{observation}

\subsection {The vector bundle $B$ on $C^{norm}$} 
Let again $d\geq 3$ denote the degree of the ruled surface $S$ and let $C\subset Gr$ be the 
corresponding curve.  Put  
\[B:=\{(p,v)\in C^{norm} \times V|\ p\mapsto \overline{w}\in C,  w\wedge v=0\}
\subset C^{norm} \times V.\]  
This is a (geometric) vector bundle of rank two on $C^{norm}$. We will identify $B$ with its sheaf 
of sections. We note that  $Proj(B)=\tilde{\tilde{S}}$. The line bundle $\Lambda ^2B$ on $C^{norm}$ 
is the pullback of the restriction of $O_{\mathbb{P}(\Lambda ^2V)}(-1)$ to $C$ and has therefore 
degree $-d$. The vector space $H^0(C^{norm},B)=0$, otherwise all the lines of $C$ pass 
through one point and $S$ is a cone. \\

{\it The vector bundle $B$ is an important tool in case $C^{norm}$ has genus 0}. \\ 
{\it For the case $d=3$} it is easily seen that $C^{norm}$ has genus 0. Let $t$ parametrize 
$C^{norm}$. Then $B$ is isomorphic to $O_{C^{norm}}(-1)\oplus O_{C^{norm}}(-2)$. 
In particular, $\tilde{\tilde{S}}$ is isomorphic to $\mathbb{P}^2$ with one point blown up 
(see \cite{Har}, V, Example 2.11.5). The sections of $B$ with a pole of order 1 at 
$t=\infty$ are $\mathbb{C}a$ and those with a pole of order $\leq 2$ at $t=\infty$ are 
$\mathbb{C}a+\mathbb{C}b$. By choosing a suitable basis of $V$ one can normalize to 
the following two cases:\\
$a=(1,t,0,0),\ b=(0,0,1,t^2)$ and $S$ has the equation $t_1^2t_4-t_2^2t_3=0$; \\
$a=(1,t,0,0),\ b=(0,1,t,t^2)$ and $S$ has  the equation $t_3^3+t_4(t_1t_4-t_2t_3)=0$.\\
This gives the classification of the ruled cubic surfaces over, say, $\mathbb{C}$. 
In Section~\ref{section2} we will follow another  method to obtain the classification of ruled cubic
surfaces over any field and compare this with Dolgachev's  method.\\
{\it For  $d=4$ and assuming that $C^{norm}$ has genus 0}, there are two possibilities 
for the  vector bundle $B$, namely:\\ 
$B\cong O_{C^{norm}}(-1)\oplus O_{C^{norm}}(-3)$ and $\tilde{\tilde{S}}$ is the Hirzebruch surface 
$\Sigma _2$,\\  
$B\cong O_{C^{norm}}(-2)\oplus O_{C^{norm}}(-2)$ and $\tilde{\tilde{S}}$ is  
$\mathbb{P}^1\times \mathbb{P}^1$.\\ 
We note in passing that the first possibility was overlooked by Cremona. 
The method of Cayley can be interpreted as taking three sections of the vector bundle  
$B(d)$ for a certain values of $d\geq 1$.     Normalizing  sections of $B$ with poles of 
order 1,2,3 at $t=\infty$, by a choice of the basis of $V$ and possibly changing $t$, 
we will arrive in Subsection~\ref{3.1} at a classification of the corresponding 
ruled quartic  surfaces.

If the genus of $C^{norm}$ is 1, the vector bundle $B$ is not helpful   for the computation. 
However, $B$ and    also $\tilde{\tilde{S}}=Proj(B)$ will be identified.      

\subsection{The possibilities for the singular locus}\label{1.4}
It is helpful for the classification of the ruled surfaces to consider   $Q:=S\cap H$ with  
$H\subset \mathbb{P}(V)$ a general plane. By Bertini's theorem, $Q$ is an irreducible 
reduced curve of degree $d$.  The morphism $C\rightarrow Q$, given by 
$\overline{w}\in C\mapsto \overline{w}\cap H\in Q$, is birational.  Thus $C^{norm}$ is the 
normalization of $Q$.  The singular locus of $S$ is written  as a union of its irreducible   
components $C_i,\ i=1,\dots ,s$ of degree $d_i$ and  generic multiplicity $m_i\geq 2$. 
The curve $Q$   meets every $C_i$  with multiplicity in $d_i$ points with  multiplicity $m_i$.   
 For every singular point $q\in Q$ one defines a number $\delta _q$ which is the sum of the 
integers    $\frac{k(k-1)}{2}$ taken over the multiplicities $k$ of $q$ itself and of all the 
singular points that    occur in the successive blow ups of $q$.  
The Pl\"ucker formula  states that the genus of   the normalization $C^{norm}$ of $Q$ is  
$\frac{(d-1)(d-2)}{2}-\sum \delta _q$.          

 For $d=3$, there is a single singular point $q$ and $\delta _q=1$ (and $q$ is a node or  cusp).      
The singular locus is described by $s=1, d_1=1, m_1=2$.            

For $d=4$, there are more possibilities.    The  singularities of a simple plane curve 
(i.e., reduced, multiplicity $\leq 3$ and in the blow ups   
 there are only singularities of multiplicity $\leq 3$)  are classified, see \cite{Bar}, p. 62,    
by formal standard equations $F\in K[[x,y]]$.  The condition that $Q$ is irreducible, has degree 4 
and    the genus of its normalization $C^{norm}$ is 0 or 1, leads to    
the list of possibilities (with their symbols or names):\\   
for $m=2$: 
\[A_2:\   x^2-y^2,\delta =1;\ \  A_3:\  x^2-y^3, \delta =1;\ \  
A_4:\  x^2-y^4,\delta =2;\ \ \ \ \ \ \ \ \ \ \  \]   
for $m=3$: 
\[D_4:\  y(x^2-y^2), \delta =3;\ \ D_5:\  y(x^2-y^3),\delta =3; \ \ 
E_6:\   x^3-y^4,\delta =3 \]    
\[ \mbox{ and the last case }E_7:\  x(x^2-y^3), \mbox{ which is ruled out by } \delta =4.
\ \ \ \ \ \ \  \ \ \ \ \ \ \ \ \ \ \ \ \ \ \ \ \ \ \ \  \  \]   
The inequalities  $\sum _{i=1}^sd_i\frac{m_i(m_i-1)}{2}\leq \sum \delta _q\leq \frac{(4-1)(4-2)}{2}$  
lead to a     list of possibilities for the singular locus, again with Cayley's convention that $d^m$
stands for an irreducible curve of degree $d$ and with multiplicity $m$ and `int' meaning
intersecting lines: 
$ 1^2$; $\ 1^3$;  $\ 2^2$; $\ 3^2$;  $\ 1^2,1^2\  \mbox{int}$; $\  1^2,1^2$;   
$\ 1^2,2^2$;  $\ 1^2,1^2,1^2$.     

\subsection{The reciprocal of a ruled surface}\label{recip}       
As before, $V$ is a vector space of dimension 4 over a field $K$.  One identifies 
$\Lambda ^4V$ with $K$.  The nondegenerate symmetric bilinear form on  
$\Lambda ^2V$, given by $(w_1,w_2)=w_1\wedge w_2\in \Lambda ^4V= K$, yields an isomorphism 
$f:\Lambda ^2V \rightarrow \Lambda ^2V^*=(\Lambda ^2V)^*$ by $f(w_1)(w_2)=w_1\wedge w_2\in K$. 
This isomorphism maps decomposable vectors of $\Lambda ^2V$ to decomposable vectors of 
$\Lambda ^2V^*$.   

Indeed, consider $f(v_1\wedge v_2)$.  Let $v_1,v_2,v_3,v_4$ be a basis of $V$. 
The kernel of $f(v_1\wedge v_2)$ has basis 
$v_1\wedge v_2,\ v_1\wedge v_3,\ v_1\wedge v_4,\ v_2\wedge v_3,\ v_2\wedge v_4$.
Let $\ell _1,\ell _2$ be a basis of $(V/Kv_1+Kv_2)^*\subset V^*$. Then $\ell _1\wedge \ell _2$ has
the same basis vectors in the kernel. Hence $f(v_1\wedge v_2)$ is a multiple of 
$\ell _1\wedge \ell _2$.

Thus $f$ induces an isomorphism $\tilde{f}:Gr(2,V)\rightarrow Gr(2,V^*)$. Let $P\subset V$ be a 
2-dimensional subspace. Then $\tilde{f}(P)$ is the 2-dimensional subspace $(V/P)^*$ of $V^*$.  
For 2-dimensional subspaces $P_1,P_2\subset V$ with $P_1\cap P_2\neq 0$ one has
$(V/P_1)^*\cap (V/P_2)^*\neq 0$. This also follows from the formula
$f(w_1)\wedge f(w_2)=w_1\wedge w_2$ for any $w_1,w_2\in \Lambda ^2V$
(for a suitable identification of $\Lambda ^4V^*$ with $K$).
 
Any 1-dimensional subspace $L\subset V$ determines the plane in $Gr(2,V^*)$
consisting of all 2-dimensional $P\subset V,\ P\supset L$ (an $\omega$-plane in \cite{Ed}). 
The image of this plane  under $\tilde{f}$ is the plane in  $Gr(2,V^*)$ consisting of all 
2-dimensional  $Q\subset (V/L)^*\subset V^*$. Since the latter is a plane of the `opposite type' 
(a $\rho$-plane in \cite{Ed}), there is no isomorphism $V\rightarrow V^*$ inducing $\tilde{f}$.  \\

Let $e_1,\dots ,e_4$ denote a basis of $V$ and $e_1^*,\dots ,e_4^*$ the dual basis of $V^*$. 
Then $e_{ij}:=e_i\wedge e_j,\ i<j$ is a basis of $\Lambda ^2V$ and $e^*_{ij}=e_i^*\wedge e_j^*,\ i<j$ 
is a basis of $\Lambda ^2V^*$. From the Pl\"ucker coordinates $\{p_{ij}\}$ of a line 
$\overline{w}\in Gr(2,V)$ one easily finds the Pl\"ucker coordinates of $f(\overline{w})\in Gr(2,V^*)$ 
by the  identities
 
\[f(e_{12})=e^*_{34}; f(e_{13})=-e^*_{24}; f(e_{14})=e^*_{23}; f(e_{23})=e^*_{14};\]
\[ f(e_{24})=-e^*_{13}; 
f(e_{34})=e^*_{12} .\]

Let an irreducible ruled surface $S\subset \mathbb{P}(V)$ (of degree $d$) be given by an 
irreducible curve $C\subset Gr(2,V)$ of degree $d$. Consider a nonsingular point 
$s\in S$ lying on a single line $\ell $ of the surface. The tangent plane  $T_{S,s}$ contains the 
line $\ell$. The same holds for the tangent planes $T_{S,s'}$ for almost all points 
$ s'\in \ell$. Since $T_{S,s'}$ cannot be all equal,  the reciprocal (or dual) surface contains all 
planes $H\supset \ell$. It now follows that the reciprocal surface 
$\check{S}\subset \mathbb{P}(V^*)$ is ruled. The corresponding curve in $Gr(2,V^*)$ is $f(C)$. 
It has also degree $d$ since the degree of the curve $f(C)$ is $d$. 
  
Using  Pl\"ucker coordinates one easily finds $\check{S}$. 
Another useful 
{\it computation of the reciprocal surface} is the following.  
Consider 
$\tilde{S}=\{(\overline{w},\overline{v})\in C\times \mathbb{P}(V)|\ v\wedge v=0\}\rightarrow 
S\subset \mathbb{P}(V)$ and a nonsingular point $\overline{v}_0\in S$ and the nonsingular point 
$\overline{w}_0\in C$ with $w_0\wedge v_0=0$. 
(We note that the tangent plane of $\tilde{S}$ at the point $(\overline{w}_0,\overline{v}_0)$ is 
mapped isomorphically to the tangent plane of $S$ at the point $\overline{v}_0$. 
The first tangent plane is the product of the tangent line of $C$ at the point $\overline{w}_0$ 
and the line $\overline{w}_0$).

Let a (local) parametrization $t\mapsto w(t)$ for $C$ be given, such that $w_0=w(t_0)$. 
Choose a decomposition $w(t)=a(t)\wedge b(t)$, locally at $t_0$. Then 
$v_0=s_0a(t_0)+(1-s_0)b(t_0)$ and $S$ has the local parametrization
$(t,s)\mapsto sa(t)+(1-s)b(t)$. The linearization of this, i.e.,
\[v_0+(s-s_0)(a(t_0)-b(t_0))+(t-t_0)(s_0a'(t_0)+(1-s_0)b'(t_0)),\] 
is a parametrization of the tangent plane $T_{S,\overline{v}_0}$.  This corresponds with the 
3-dimensional subspace of $V$ with basis 
\[v_0, v_0+a(t_0)-b(t_0),v_0+s_0a'(t_0)+(1-s_0)b'(t_0).\]
The exterior product 
$a(t_0)\wedge b(t_0)\wedge (s_0a'(t_0)+(1-s_0)b'(t_0))$ of these vectors is an element of 
$\Lambda ^3V=V^*$. This defines a point in $\mathbb{P}(V^*)$ corresponding to the tangent plane 
$T_{S,\overline{v}_0}$. The reciprocal surface $\check{S}$ consists of all these points. 
In varying $s_0$ one finds a line on $\check{S}$, through the points
$a(t_0)\wedge b(t_0)\wedge a'(t_0)$ and $a(t_0)\wedge b(t_0)\wedge b'(t_0)$.
By varying $t_0$ one obtains $\check{S}$.

For the two cases of ruled quartic surfaces $S$ with genus 1, it is easily seen that 
$\check{S}\cong S$.  For the ruled quartic surfaces cases of genus 0  there are explicit global
 expressions $w(t)=a(t)\wedge b(t)$ (with $t\in \mathbb{P}^1$) and the above can be used for the 
computation of $\check{S}$.\\     

\subsection{The classification of quartic ruled surfaces}\label{s1.6}     
The {\it Number} appearing in the tables are introduced for notational reasons in the computation 
of Subsection~\ref{3.1}. It has no intrinsic meaning.   
The cases with $C^{norm}$ of genus 0 and $B$ of type $-1,-3$ are:\\  
      
\noindent     $ \begin{array}{|c|c|c|c|c|c|c|}  \hline 
 \mbox{Number}& \mbox{singular}& \dim & \mbox{singularities }& \mbox{tangent}&
\mbox{Cremona}& {\rm XIII}\\    
& \mbox{locus}& P(C) & \mbox{of }C& \mbox{spaces}& & \\    \hline      
1\ a,b,c & 1^3&3& - & 2 & 9&\\   \hline    
2\ a,b,c & 3^2&4&-&1&7&\\    \hline    
3\ a,b  & 1^2,2^2&4&-&1&4& \\    \hline    
4& 1^3&3& \mbox{node}& 1& 10&7 \\     \hline   
 5&1^3&3&\mbox{cusp}&1&10&\\    \hline    
\end{array} $\\            

The cases with $C^{norm}$ of genus 0 and $B$ of type $-2,-2$ are:\\  
      
\noindent    $ \begin{array}{|c|c|c|c|c|c|c|}  \hline  
\mbox{Number}& \mbox{singular}& \dim & \mbox{singularities}& \mbox{tangent }&
\mbox{Cremona}& {\rm XIII}\\      
& \mbox{locus}& P(C) & \mbox{of }C& \mbox{spaces}& &\\    \hline      
6&1^2,1^2,1^2& 3&\mbox{node} &2&5&\\   \hline  
 7&1^2,1^2,1^2& 3&\mbox{cusp}&2&5&\\   \hline   
8& 1^3& 4& - & 1& 3&\\   \hline   
9& 1^3&4&-&1&3&\\   \hline   
10\ a & 1^2,1^2, int & 3 &\mbox{node}& 1&6&\\   \hline   
10\  b & 1^2,1^2, int & 3 &\mbox{ cusp}& 1&6&\\   \hline   
11 & 1^2,2^2 & 4 & - & 1& 2&\\   \hline     
12 & 1^2,2^2 & 4 & - & 1& 2&8 \\    \hline    
13 \ a,b,c & 1^3&4&-& 1& 8&6 \\    \hline   
 14 \ ... & 3^2& 4 & - & - & 1&9,10\\    \hline   
 \end{array} $ \\
 
In \S~\ref{3.2.4}  it is shown that {\it Number} 14 consists of six distinct cases.  
In the text it is explained that the cases 11 and 12 are different. The reciprocals of
$2\ a,b,c$ are $13\ a,b,c$ and the reciprocals of  $3\ a,b$ are $8, 9$. The other examples are 
`selfdual'. 

The cases with $C^{norm}$ of genus 1 are:\\

\noindent$ \begin{array}{|c|c|c|c|c|c|}  \hline  
\mbox{Number}& \mbox{ singular locus} & \dim P(C)& \mbox{tangent spaces}&\mbox{Cremona}&
{\rm XIII}\\   \hline  
15& 1^2& 3 & 1&12& 5\\       \hline  
16 & 1^2,1^2 & 3 &  2 & 11& 1,2,3,4 \\      
\hline    \end{array} $    

 \section{Ruled surfaces of degree 3}\label{section2}
Here we give the classification over an arbitrary field $K$. The singular locus of $S$ is a line,  
$C^{norm}$ has genus 0 and $\dim P(C)=3$.  This implies that $C=C^{norm}$ is the twisted cubic 
curve in $P(C)$.  \\

{\it In the first case $P(C)$ lies in two tangent spaces} at the points 
$\overline{w}_1,\overline{w}_2\in Gr$. From  Lemma~\ref{1.2.10} we conclude that
$S$ is given by a bi-homogeneous equation $F$ in the pairs of variables 
$t_1,t_2$ and $t_3,t_4$ of bi-degree $(2,1)$, corresponding to a morphism 
$f:\overline{w}_1\rightarrow \overline{w}_2$ of degree 2. The line $\overline{w}_1$ is nonsingular 
and a `directrix'. The line $\overline{w}_2$ is the singular locus. Further
$\tilde{C}\setminus C=\{\overline{w}_1,\overline{w}_2\}$.

If the field $K$ has characteristic $\neq 2$, then $C,P(C),\overline{w}_1,\overline{w}_2$ are all 
defined over $K$ and can be put in a standard form. The morphism $f$ is defined over $K$. 
The ramification points of $f$ are either both in $K$ or belong to a quadratic extension of $K$ 
and are conjugated.    

If the field $K$ has characteristic 2, then $f$ is either separable and has one point of ramification, 
or $f$ is inseparable. In both cases $f$ can be put into a standard form.\\    
  
{\it In the second case $P(C)$ lies in only one tangent space}, namely at the point 
$\overline{w}_0$ which is   the singular line of $S$. Then $C$ lies on the quadratic cone in 
$P(C)$ and $\overline{w}_0\in C$. In this   case $\tilde{C}=C$. Now $C$ and $S$ can be put into a 
standard form. We arrive at the following result.           

\begin{proposition}\label{cubic}\label{2.0.1} The standard equations for ruled cubic surfaces $S/K$, which 
are not cones, are the following:\\  
{\rm (1)}  $t_3t_1^2+t_4t_2^2=0$. If $char K=2$, then there are no twist. 
For $char K\neq 2$  the twists are $t_3(t_1t_2)+t_4(at_1^2+t_2^2)=0$ with $a\in K^*$ not a square.\\  
{\rm (2)} $t_3t_1t_2+t_4t_1^2+t_2^3=0$ (there are no twists).\\  
{\rm (3)} $t_3t_1^2+t_4(t_2^2+t_1t_2)=0$ if  $char K=2$ (there are no twists).  
\end{proposition} 

The curves  $C$ for (1) and (2)  are in parameter form 
\[p_{12}=0,\ p_{13}=-t^2,\ p_{14}=1,\ p_{23}=-t^3,\  p_{24}=t,\ p_{34}=0
\mbox{ and }\]
\[ p_{12}=0,\ p_{13}=t^3,\ p_{14}=t^2,\ p_{23}=-t^2,\  p_{24}=-t,\ p_{34}=-1.\]
The above equations for $S$ are derived in an elegant way by I.~Dolgachev \cite{Dol}, using only 
the information that the singular locus of $S$ is a line with multiplicity 2.

For $K=\mathbb{R}$, there are three types of cubic ruled surfaces (omitting  cones). 
Case (1) of Proposition~\ref{cubic}, without twist, is represented by the plaster model VII, nr. 21 
and by the string models XVIII, nr.~2 and 3.
The twisted case ($a=-1$) is represented by VII, nr.~20 and XVIII, nr.~1.
Finally, case (2) carries the name Cayley's ruled cubic surface. It is
represented by VII, nr.~22 and 23 and XVIII, nr.~4.                

\section{Ruled surfaces of degree 4}\label{section3}
The base field $K$ is supposed to be algebraically closed. The only role
that the characteristic of $K$ plays is in the classification of the morphisms 
$\mathbb{P}^1\rightarrow \mathbb{P}^1$ of degree 2 and 3. For convenience we suppose that 
$K$ has characteristic 0 or $>3$. We need both the vector bundle $B$ and the possibilities 
for the singular locus in order to find all cases and to verify the computations by comparison.  

\subsection{Classification of $S$, 
using the vector bundle $B$}\label{3.1}
\subsubsection{$C^{norm}$ of genus 0 and $B$ of type $-1,-3$}\label{3.1.1}
Choose a parameter $t$ for $C^{norm}\cong \mathbb{P}^1$ and let $p$ be the pole of
$t$. Then $H^0(C^{norm}, B([p])$ has basis $a$ and $H^0(C^{norm},B(3\cdot [p]))$ has basis 
$a,ta,t^2a,b$. Now $a,b\in K[t]\otimes _KV$  have degrees 1 and 3.
We note that $b$ is not unique and can be replaced by 
$\mu b+\lambda _0a+\lambda _1ta+\lambda _2t^2a$ with $\mu \neq 0$.
We will derive {\it normal forms} for the possibilities of the pair $a,b$. These will depend
on the choice of $t$. 
There is a unique subspace $W\subset V$ of dimension 2 with
$a\in K[t]\otimes W$ and the image $b'$ of $b$ in $K[t]\otimes V/W$ is unique up to
multiplication by some $\mu \in K^*$ and has degree $\geq 1$.

The above normal form  is obtained by choosing $t$ and  a basis $e_1,e_2,e_3,e_4$ of $V$ such 
that $\{e_1,e_2\}$ is a basis of $W$ and such that $a=(a_1,a_2,0,0)$ and a suitable 
$b=(b_1,b_2,b_3,b_4)$ w.r.t. this basis are as simple as possible. 

The Pl\"ucker coordinates of the line through $a(t)$ and $b(t)$ are easily computed and this yields  
$C\subset Gr$ in parametrized form. From this one deduces $\dim P(C)$, possible singularities 
of $C$ and the relation of $C$ w.r.t. the tangent spaces of $Gr$. The reciprocal surface 
(needed for the comparison with Cremona's list) is computed by the method of 
Subsection~\ref{recip}, 
again using $a(t),b(t)$.  For some cases (especially when the singular locus is $3^2$), 
the equation of the surface $S$ is rather long and requires a MAPLE computation. 
We avoid this and compute the singular locus by other means. 
We start by classifying the pairs $(b_3,b_4)$ which are uniquely determined by $S$, up to taking 
linear combinations. \\

\noindent  $\max (\deg b_3, \deg b_4)=3$ and $\gcd(b_3,b_4)=1$.\\
The morphism $C^{norm}\rightarrow \mathbb{P}^1, t\mapsto (b_3:b_4)$ has degree 3. 
The possibilities for the ramification indices are: (a) $3,3$, (b)  $3,2,2$ or (c) $2,2,2,2$. 
A  change of $t$ and a linear change of $e_3,e_4$ brings the pair $(b_3,b_4)$ into a standard form
\[(1,t^3),\  (1,t^2(t+1)), \mbox{ or } ( t-\mu , (2\mu -1)t^3+(2-3\mu )t^2 ) \mbox{ with } \mu \neq 0,1,1/2.\]  
(In the last case the 4 ramification points are $0,1,\infty ,\frac{\mu}{2\mu -1}$).
One normalizes $a=(1,t,0,0)$,  $b=(b_1,b_2,b_3,b_4)$ and $\max (\deg b_1,\deg b_2)\leq 2$. 
Then $b-b_1\cdot a$ produces a new $b=(0,b_2,b_3,b_4)$. There are now again  two cases:

\paragraph{1.} (a,b,c).  $b_2=0$ and the data are: $C$ is nonsingular,  
$P(C)=T_{Gr,\overline{e_{12}}}\cap T_{Gr,\overline{e_{34}}}$  and
$\overline{e_{12}},\overline{e_{34}}\not \in C$; 
equation $t_1^3b_3(\frac{t_2}{t_1})t_4-t_1^3b_4(\frac{t_2}{t_1})t_3=0$; the singular locus of $S$ is 
the line $\overline{e_{34}}$ with multiplicity 3. Then $1^3$ and Cremona 9.

\paragraph{2.} (a,b,c). $b_2\neq 0$ and the data are: $C$ nonsingular, 
$P(C)=T_{Gr,\overline{e_{12}}}$ and $\overline{e_{12}}\not \in C$.  A direct computation of the 
equation seems difficult. The points of the line through $a(t)$ and $b(t)$ can be written as
$(1,t+\lambda \frac{b_2}{b_3}(t),\lambda ,\lambda \frac{b_4}{b_3}(t))$. Computing with the equality
\[  (1,t+\lambda \frac{b_2}{b_3}(t),\lambda ,\lambda \frac{b_4}{b_3}(t))=
(1,s+\lambda \frac{b_2}{b_3}(s),\lambda ,\lambda \frac{b_4}{b_3}(s))\mbox{, with } s\neq t\]  
leads to the result that the `double  curve' and thus also the singular locus is a twisted cubic curve.  
Then $3^2$ and Cremona 7.

\paragraph{3.}  (a,b) $\max (\deg b_3,\deg b_4)= 2$ and $\gcd(b_3,b_4)=1$. \\
A normalization is $a=(t+\beta ,1,0,0)$ and $b=(0,t^3+\alpha t,t^2,1)$. Equation 
\[ t_3t_4(t_2-\beta (t_3+\alpha t_4))^2-(t_3(t_3+\alpha t_4)-\beta t_2t_4+t_1t_4)^2=0.\]
One has $P(C)=T_{Gr,\overline{e_{12}}}$ and $\overline{e_{12}}\in C$.
The singular locus is the union of the line $\overline{e_{12}}$ and the conic
$t_2-\beta (t_3+\alpha t_4)=0, \ t_3(t_3+\alpha t_4)-\beta t_2t_4+t_1t_4=0$.
Then  $1^2,2^2$ and Cremona 4. The degree morphism $C\rightarrow D$ has
two points of ramification. The point $L\cap D$ is a ramification point on $D$ if and only if 
$\beta =0$. We consider this as two cases.

\paragraph{4.} $\max (\deg b_3,\deg b_4)= 2$ and $\gcd(b_3,b_4)$ has degree 1.\\
A normalization  is  $a=(t,1,0,0),\ b=(0,t^3+\alpha ,t(t+\beta),t+\beta )$
with $\gcd(t^3+\alpha ,t+\beta )=1$. The equation is
\[t_1t_4^2(t_3+\beta t_4) -t_2t_3t_4(t_3+\beta  t_4)+\alpha t_3t_4^3+t_3^4=0 .\]
Further $\overline{e_{12}}\in C$ is a node (for $t=\infty , t=-\beta$), $\dim P(C)=3$ and $P(C)$ 
lies in only one tangent space $T_{Gr,\overline{e_{12}}}$. Then $1^3$ and Cremona 10.

\paragraph{5.} $\max (\deg b_3,\deg b_4)= 1$.\\
A normalization is $a=(t,1,0,0),\ b=(0,t^3+\alpha t^2,t,1)$. The equation is
\[t_1t_4^3-t_2t_3t_4^2+\alpha t_3^3t_4+t_3^4=0.\]
Further $\overline{e_{12}}\in C$ is a cusp (for $t=\infty$), $\dim P(C)=3$ and $P(C)$ lies in only 
one tangent space, namely $T_{Gr,\overline{e_{12}}}$. Then $1^3$ and Cremona 10.\\

Finally we show that the omitted cases  can be reduced to the above.\\ 

$\max (\deg b_3,\deg b_4)=3$ and $\gcd(b_3,b_4)$ has degree 1.\\
A normalization is $a=(1,t,0,0),\ b=(b_1,b_2,t, t(t+\mu )^2)$. Replacing $t$ by 
$s^{-1}$ and 
multiplying by a suitable power of $s$ yields
$a=(s,1,0,0),\ b=(s^3b_1(s^{-1}),s^3b_2(s^{-1}), s^2, (1+\mu s)^2)$. 
Thus reduction to $\max (\deg b_3,\deg b_4)= 2$.\\
 
$\max (\deg b_3,\deg b_4)=3$ and $\gcd(b_3,b_4)$ has degree 2.\\
A normalization is $a=(1,t,0,0),\ b=(b_1,b_2,t(t+\mu),t(t+\mu)(t+\lambda ))$.
Replacing $t$ by $s^{-1}$ and multiplying by a suitable power of $s$ gives a reduction to  
$\max (\deg b_3,\deg b_4)= 2$.

\subsubsection{$C^{norm}$ of genus 0 and $B$ of type $-2,-2$}\label{3.1.2}
$V,t,p$  have the same meaning as in \S~\ref{3.1.1}. The vector space $H^0(C^{norm},B(2[p]))$ has 
dimension 2 and consists of elements in $K[t]\otimes V$ of degree $\leq 2$ and the only 
element of degree $\leq 1$ is 0.
We are interested in lines $Ka\subset H^0(C^{norm},B(2[p]))$ such that 
$a\in K[t]\otimes W$ with $\dim W=2$.\\

\noindent {\it Suppose that there are two such lines $Ka$ and $Kb$}. \\

One can normalize by $a=(a_1,a_2,0,0),\ b=(0,0,b_3,b_4)$. The two morphisms
$\overline{a}, \overline{b}:C^{norm}\rightarrow \mathbb{P}^1$, $t\mapsto (a_1:a_2)$ and 
$t\mapsto (b_3,b_4)$ of degree 2 are distinct and their sets of ramification points can be either 
disjoint or have one point of intersection. Choosing $t$ leads to the following normalizations.

\paragraph{6.}   $a=(1,t^2,0,0), \ b=(0,0,(t-1)^2,(t-\lambda ^2 )^2)$. The singular locus consists 
of the lines $\overline{e_{12}},\overline{e_{34}}$ and a third line corresponding to 
$t=\pm \lambda$. The morphism 
$C^{norm}\stackrel{(\overline{a},\overline{b})}{\rightarrow} \mathbb{P}^1\times \mathbb{P}^1$  
maps $t=\pm \lambda$ to the same point of $C$. Thus $C$ has a node, 
$P(C)=T_{Gr,\overline{e_{12}}}\cap T_{Gr,\overline{e_{34}}}$,  $1^2,1^2,1^2$ and Cremona 5.

\paragraph{7.}  $a=(1,t^2,0,0),\ b=(0,0,1,(t-1)^2)$. The image $C$ of 
$C^{norm}\stackrel{(\overline{a},\overline{b})}{\rightarrow} \mathbb{P}^1\times \mathbb{P}^1$  
has a cusp corresponding to  $t=\infty$. The singular locus consists of three lines 
$\overline{e_{12}},\overline{e_{34}}$ and the line corresponding to $t=\infty$.
Thus  $P(C)=T_{Gr,\overline{e_{12}}}\cap  T_{Gr,\overline{e_{34}}}$,  $1^2,1^2,1^2$ and 
Cremona 5.\\

\noindent {\it Suppose that there exists only one such line $Ka$}. \\

Normalize by $a=(1,t^2,0,0),\ b=(b_1,b_2,b_3,b_4)$ with $\deg b_2<2$. The pair $(b_3,b_4)$ is, 
up to taking linear combinations, uniquely determined by the surface.
The morphism  $m: C^{norm}\rightarrow \mathbb{P}^1,\ t\mapsto (b_3:b_4)$  cannot be constant 
and has degree 1 or 2. There are the following cases.  

\paragraph{8.}  $(b_3,b_4)=(1,t)$. Then $P(C)=T_{Gr,\overline{e_{12}}}$ and  
$\overline{e_{12}}\not \in C$. The equation is
\[t_1t_3t_4^2-t_2t_3^3-t_3^2t_4^2b_1(\frac{t_4}{t_3})+t_3^4b_2(\frac{t_4}{t_3})=0.\] 
Thus $1^3$ and Cremona 3.

\paragraph{9.}  $(b_3,b_4)=(t-\alpha ,t(t-\alpha ))$. Then 
$P(C)=T_{Gr,\overline{e_{12}}}$ and $\overline{e_{12}} \in C$. The equation is 
\[t_1(t_4-\alpha t_3)t_4^2-t_2t_3^2(t_4-\alpha t_3)-t_3^2t_4^2b_1(\frac{t_4}{t_3})+
t_3^4b_2(\frac{t_4}{t_3})=0.\]
Thus $1^3$ and Cremona 3. 

\paragraph{10.} (a,b)  Now the morphism $m$ has degree two. If $t=0,\infty$ are the ramification 
points of $m$, then one normalizes to $(b_3,b_4)=(1,t^2)$. Then $\dim P(C)=3$ and $P(C)$ lies in 
only one tangent space, namely $T_{Gr,\overline{e_{12}}}$, and $\overline{e_{12}}\not \in C$.  
Write  $b_1=b_{12}t^2+b_{11}t+b_{10}$ and $b_2=b_{21}t+b_{20}$. 
One can normalize further to $b_1=b_{11}t,\ b_2=b_{21}t$. Then $C$ has a node (case (a)) if 
$b_{21}b_{11}\neq 0$ and has a cusp otherwise (case (b)).  Then $1^2,1^2,\ int$ and Cremona 6. 

\paragraph{11.}  If $m:C^{norm}\rightarrow \mathbb{P}^1$ is ramified for, say, $t=1,\infty$, then one 
can normalize $a=(1,t^2,0,0),\ b=(b_1(t-1),b_2(t-1),1,(t-1)^2))$ with $b_1,b_2\in K$. The equation is
\[t_3t_4(2t_1+(b_2-b_1)t_3-b_1t_4)^2-(t_2t_3-t_1t_3-t_1t_4+2b_1t_3t_4)^2=0.\]
The singular locus is the union of the line $L=\overline{e_{12}}$ and the conic
$D$ given by $2t_1+(b_2-b_1)t_3-b_1t_4=0, \ t_2t_3-t_1t_3-t_1t_4+2b_1t_3t_4=0$.
Now $P(C)=T_{Gr,\overline{e_{12}}}$, $\overline{e_{12}}\not \in C$, the image of 
$C\rightarrow L\times D,\ \overline{w}\mapsto (\overline{w}\cap L,\overline{w}\cap D)$ is a rational 
curve having a cusp. Then $1^2,2^2$ and Cremona 2.

\paragraph{12.}  If $m$ is ramified for, say, $t=1, \mu$ with $\mu \neq 0,1,\infty $, then one can 
normalize to $a=(1,t^2,0,0)\ b=(b_1,b_2,(t-1)^2,(t-\mu)^2)$ with $b_1,b_2\in K$.\\
One can replace $b$ by $b-b_1\cdot a$ and normalize further to
$a=(1,t^2,0,0),\ b=(0,1,(t-1)^2,(t-\mu )^2)$. A somewhat long computation yields the equation
\[4t_3t_4( (\mu-1)^2(t_2-\mu t_1)-2t_3-2t_4)^2-\]
\[((\mu -1)^2(-\mu ^2 t_1t_3-t_1t_4+t_2t_3+t_2t_4)-t_3^2-6t_3t_4-t_4^2)^2=0.\]
The singular locus is the union of the line $L=\overline{e_{12}}$ and the conic $D$
given by the equations 
\[(\mu-1)^2(t_2-\mu t_1)-2(t_3+t_4)=0, \ (\mu-1)^3t_1(t_3-\mu t_4)-(t_3-t_4)^2=0.\]
Further $P(C)=T_{Gr,\overline{e_{12}}}$, $\overline{e_{12}}\not \in C$,
the image of the morphism $C\rightarrow L\times D$ is a rational curve  having a node. 
Then $1^2,2^2$ and Cremona 2. \\

\noindent  {\it  Suppose that  there is no such line and that 
$P(C)$ lies in a tangent space}.\\

The inclusion $P(C)\subset T_{Gr,\overline{e_{12}}}$ yields a morphism 
$f:C^{norm}\rightarrow \overline{e_{12}}$  induced by  
$\overline{e_{12}}\neq \overline{w}\in C  \mapsto \overline{w}\cap \overline{e_{12}}$.  
If the degree of $f$ is $1$, then we may assume that $(1,t,0,0)$ lies on $S$. Combining with a 
nonzero element  $a\in H^0(C^{norm},B(2[p]))$, one finds a surface of degree 3 instead of 4. \\

\noindent The possibility that {\it the degree of $f$ is $2$} is excluded by the following  reasoning. 
Let $t$ be a parameter for $C^{norm}$ and write $f=((\alpha t+\beta )^2,(\gamma t+\delta )^2,0,0)$. 
Let  $a,b$ be  a basis of $H^0(C^{norm},B(2[p]))$. Then 
$\lambda _0(t)f=\lambda _1(t)a+\lambda _2(t)b$ holds for some 
$\lambda _0(t),\lambda _1(t),\lambda _2(t)\in K[t]$ with $\gcd(\lambda _1(t),\lambda _2(t))=1$.
The Pl\"ucker coordinates of $a\wedge b$ are polynomials in $t$ with 
greatest common divisor 1 and
maximal degree 4, since these parametrize $C$. The same holds for the Pl\"ucker coordinates of 
$f\wedge a$ and $f\wedge b$. The equality 
$\lambda _0(t)f\wedge a=-\lambda _2(t)\cdot  a\wedge b$ implies that $\lambda _0(t)$ is a constant 
multiple of $\lambda _2(t)$. Similarly, $\lambda _0(t)$ is a constant multiple of $\lambda _1(t)$.
We conclude that the $\lambda _i(t)$ are constant. Then$f\in H^0(C^{norm},B(2[p])$ and this 
contradicts the assumption.\\ 

\paragraph{13.} (a,b,c). {\it If the degree of $f$ is $3$},  then $\overline{e_{12}}$ has multiplicity 3 
and thus $1^3$. As in case 2, there are three different possibilities for the ramification of $f$. 
One writes $f(t)=(c_1,c_2,0,0)$ where $c_1,c_2$ are relatively prime polynomials in $t$ and, say, 
$\deg c_1 <\deg c_2=3$. Let $a(t)=(a_1,a_2,a_3,a_4)$ be a nonzero section of $B(2[p])$. 
An inspection of the Pl\"ucker coordinates of $f\wedge a$ implies that 
$\max (\deg a_3, \deg a_4)\leq 1$.  Moreover $a_3,a_4$ are linearly independent. 
Thus we may normalize to $(a_3,a_4)=(1,t)$.
Because $\overline{e_{12}}$ has multiplicity 3, the equation for $S$ has the form 
$t_1A_1+t_2A_2+A_3=0$, where $A_1,A_2,A_3$ are homogeneous polynomials in $t_3,t_4$ of 
degrees $3,3,4$. Substitution of  $(\lambda c_1+a_1,\lambda c_2+a_2,1,t)$ in this equation yields 
$c_1(t)A_1(1,t)+c_2(t)A_2(1,t)=0$ and we can normalize to $A_1(1,t)=c_2(t),\ A_2(1,t)=-c_1(t)$. 
In particular, $\gcd(A_1,A_2)=1$. Further
$A_3(1,t)=-a_1(t)c_2(t)+a_2(t)c_1(t)$. The term $A_3$ cannot be made 0 by a transformation of 
the form $t_1\mapsto t_1+*t_3+*t_4, \ t_2\mapsto t_2+*t_3+*t_4$, since $P(C)$ does not lie in 
another tangent space. Therefore, $\max (\deg a_1,\deg a_2)=2$.  Further $\dim P(C)=4$ and 
$\overline{e_{12}}\not \in C$. One verifies that
the equations belong to the case that $B$ has type $-2,-2$ by comparing with the cases
$1^3$ where $B$ has type $-1,-3$. Further Cremona 8. \\

\noindent {\it Suppose that there is no such line and $P(C)$ does not lie in a tangent space}.

\paragraph{14.}  {\it We claim that the singular locus is $3^2$ and is of species Cremona }1.\\
The conditions imply that $\dim P(C)=4$ and $C$ is nonsingular.  Suppose that the 
singular locus of $S$ contains a line. 
This line belongs to $C$ (because of Lemma~\ref{1.2.3}) and is, say, $\overline{w}(0)\in C$. 
Take a plane $H$ containing $\overline{w}(0)$. The intersection $H\cap S$ consists of 
$\overline{w}(0)$ with multiplicity $\geq 2$ and a remaining curve $R$ which is a conic or two 
lines or one line. For $t\neq 0$ the intersection $\overline{w}(t)\cap H$ lies on $R$. 
The possibility that $R$ is one or two lines contradicts that $P(C)$ does not lie in a tangent space. 
Thus $R$ can only be a conic. For $t\neq 0$, the positive divisor  $\overline{w}(t)\cap R$  
has degree 1 and has degree 2 for $t=0$. This is a contradiction.
  
We conclude that the singular locus of $S$ does not contain a line. 
Then, because of Lemma~\ref{3.2.1} and Subsection~\ref{1.4},  the singular locus of $S$ is  
the twisted cubic curve.  In \S~\ref{3.2.4} it is shown that this {\it Number}~14 consists of six
subclasses.

\subsubsection{The vector bundle $B$ for a genus 1 curve $C^{norm}$}\label{3.1.3}
Here we use the information from \S~\ref{3.2.5} and \S~\ref{3.2.7} below, and deduce the structure
of the vector bundle $B$ on the genus 1 curve $C=C^{norm}$.

\paragraph{15.} {\it Case $1^2$}. The equation is $(t_1t_4-t_2t_3)^2+H(t_3,t_4)$, where $H$ is 
homogeneous of degree 4  and defines 4 distinct points on 
$\mathbb{P}(Ke_3+Ke_4)=\mathbb{P}^1$. We may suppose that these points are 
$0,1,\lambda ,\infty$.  The lines $\overline{w}(t)\in C$ on $S$ are computed to be the lines 
passing through the points $(1,t,0,0)$ and $(0,y,1,t)$, with $y^2=H(1,t)$.  The genus one curve 
$C$ is made into an elliptic curve by the choice of the neutral element $e$ to correspond to 
$t=y=\infty$.  We note that $\overline{e_{12}}\not \in C$. $(1,t,0,0)$ is a section of $B(2[e])$ and 
$(0,y,1,t)$ is a section of $B(3[e])$. Further 
$w(t)=(1,t,0,0)\wedge (0,y,t,1)=ye_{12}+e_{13}+te_{14}+te_{23}+t^2e_{24}+0e_{34}$ is a section of 
$\Lambda ^2B([4]e)$.  Consider the exact sequence  
\[0\rightarrow O_C(1,t,0,0) \rightarrow B(2[e])\rightarrow O_C(0,0,1,t)\rightarrow 0.\] 
From  $O_C(1,t,0,0)\cong O_C(0,0,1,t)\cong O_C$ and $H^0(C,B(2[e]))=K(1,t,0,0)$ one concludes 
that the sequence does not split. Therefore the ruled surface (in the modern sense) 
$\tilde{S}\rightarrow C$ corresponds to the unique indecomposable vector bundle on $C$ which 
is an extension of $O_C$ by $O_C$. (see \cite{Har}).  

\paragraph{16.} {\it Case $1^2,1^2$}.  The equation $F$ for $S$ is bi-homogeneous of degree$(2,2)$ 
in the pairs of variables $t_1,t_2$ and $t_3,t_4$. The equation $F$ also defines a genus one curve
$E\subset \mathbb{P}(Ke_1+Ke_2)\times \mathbb{P}(Ke_3+Ke_4)$. Further 
$E \rightarrow C\subset Gr$ is the isomorphism which sends $p\in E$ to the line through the points 
$(pr_1(p),0,0)$ and $(0,0,pr_2(p))$. The vector bundle $B$ is the direct sum of the line bundles 
$\mathcal{L}_1:=\{ (\overline{w},v)|\overline{w}\in C,\ v\in Ke_1+Ke_2,\ w\wedge v=0\}$ and    
$\mathcal{L}_2:=\{ (\overline{w},v)|\overline{w}\in C,\ v\in Ke_3+Ke_4,\ w\wedge v=0\}$ 
of degree $-2$.  

A line bundle $\mathcal{L}$ on $E$ of degree~$-2$ induces a degree $2$ morphism
$E\rightarrow \mathbb{P}(H^0(E,\mathcal{L}^*))$. This yields a bijection between the isomorphy 
classes of line bundles of degree~$-2$ and the equivalence classes of nonconstant morphisms 
$E\rightarrow \mathbb{P}^1$ of degree~$2$. Then $\mathcal{L}_1$ is not isomorphic to 
$\mathcal{L}_2$, since the two morphisms are not equivalent.
The ruled surface $\tilde{S}\rightarrow E$ is equal to $Proj(O_E\oplus \mathcal{L})$, where 
$\mathcal{L}=\mathcal{L}_1\otimes \mathcal{L}_2^{-1}$ is any line bundle of degree 0, not 
isomorphic to $O_E$. In particular,
$\tilde{S}\not \cong \mathbb{P}^1\times E$.

\subsection{The classification, using the singular locus}\label{3.2}

\subsubsection{$2^2$ does not occur as singular locus}

\begin{lemma}\label{3.2.1} The singular locus of a quartic ruled surface cannot be  a conic. 
\end{lemma}

\begin{proof}  Suppose that the conic $D$, lying in a plane $H\subset \mathbb{P}(V)$, is the 
singular locus of some ruled quartic surface $S$, corresponding to a curve 
$C\subset \mathbb{P}(\Lambda ^2V)$. If $C$ has genus 1, then $P(C)$ is contained in a tangent 
space of $Gr$ at some point $\overline{w}_0$.
The morphism $\overline{w}\in C\mapsto \overline{w}\cap \overline{w}_0\in \overline{w}_0$ has 
degree at least 2 and thus $\overline{w}_0$ belongs to the singular locus. Hence 
$C^{norm}$ has genus 0. The morphism $f: C^{norm}\rightarrow D$, given by 
$\overline{w}\in C^{norm}\mapsto \overline{w}\cap H\in D$, has degree at  most 2, 
since the multiplicity  of $D$ is 2.\\  

Suppose that the degree of $f$ is 1. One can parametrize 
$C^{norm}$ with a parameter $t$ and choose coordinates for $\mathbb{P}(V)$ such that the line 
$\overline{w}(t)\in C^{norm}$ intersects the conic $D$ in the point $(0,1,t,t^2)$. Let $(1,0,a,b)$ with 
$a,b\in K(t)$ be another point of this line $\overline{w}(t)$. The Pl\"ucker coordinates
of $\overline{w}(t)$ are
\[p_{12}=1,\ p_{13}=t,\ p_{14}=t^2,\ p_{23}=-a,\ p_{24}=-b,\ p_{34}=t^2a-tb .\]  
Let $d$ be the common denominator of $a$ and $b$. Then $\{dp_{ij}\}$ are polynomials of degree 
$\leq 4$ and with $\gcd=1$. If $\alpha$ is a zero of $d$, then the line $\overline{w}(\alpha )$ lies in 
the plane $H$. Since this is not possible, $d=1$ and $a,b\in K[t]$. One obtains the contradiction 
that the line $\overline{w}(\infty )$ lies in the plane $H$.\\

Suppose that  the degree of $f$ is 2. One can parametrize  $C^{norm}$ with parameter $t$, and 
choose coordinates for $\mathbb{P}(V)$ such that $\overline{w}(t)\mapsto (0,1,t^2,t^4)\in D$. 
The line $\overline{w}(t)$ goes through a point $(1,0,a,b)$ where $a,b\in K(t)$. 
The Pl\"ucker coordinates of $\overline{w}(t)$ are 
\[p_{12}=1,\ p_{13}=t^2,\ p_{14}=t^4,\ p_{23}=-a,\ p_{24}=-b,\ p_{34}=t^4a-t^2b .\]
Let $d$ be the common denominator of $a$ and $b$. After multiplying the Pl\"ucker coordinates 
with $d$, the degrees are bounded by 4. Hence $d=1$ and $a,b\in K[t]$. Further the degree of $a$ 
is $\leq 2$ and the degree of $c:=b-t^2a$ is $\leq 2$. The symmetric polynomial 
$w(s)\wedge w(t)$  in $s,t$ can only have the factors $s+t$ and $s-t$. Indeed,  
$t\neq s$ and $w(s)\wedge w(t)=0$ implies that $\overline{w}(s)\cap \overline{w}(t)\in D$ and thus 
$s=-t$. It follows that $a=a_0+a_2t^2,\ c=c_0+c_2t^2$ and this contradicts that $C^{norm}$ is 
parametrized by $t$.  
\end{proof}

\subsubsection{$1^2,2^2$}\label{3.2.2}
The curve $C^{norm}$, corresponding to a ruled quartic surface $S$ with this type of singular locus, 
has genus 0 by Observation~\ref{1.2.11}. The singular locus is the union of a conic $D$ and a line $L$. 
The plane  $H\supset D$ satisfies $S\cap H=C$. Thus $L$ does not lie in $H$ and the intersection 
$L\cap H$ is a point of $D$. As in the proof of Lemma~\ref{3.2.1}, one shows that the morphism 
$C^{norm}\rightarrow D$, given by $\overline{w}\in C\mapsto \overline{w}\cap H\in D$, 
has degree 2.  Let  $D=\{(0,1,\mu^2 ,\mu )|\  \mu \in \mathbb{P}^1\}$ and
$L=\{(1 ,\lambda ,0,0)| \ \lambda \in \mathbb{P}^1\}$. The equations for $D$ and $L$
are $t_1=t_2t_3-t_4^2=0$ and $t_3=t_4=0$. The equation $F$ for $S$ lies in the ideal 
$(t_1,t_2t_3-t_4^2)^2\cap (t_3,t_4)^2$. Thus $F=t_1^2A_2+t_1(t_2t_3-t_4^2)A_1+(t_2t_3-t_4^2)^2$ 
where $A_2$ and $A_1$ are homogeneous of degrees 2 and 1.  One may suppose that $A_1$ does 
not contain $t_1$.
If $A_1$ contains $t_2$, then $F$ contains the monomial $t_1t_2^2t_3$ which is not possible. 
Hence $A_1$ is linear in $t_3,t_4$ and it follows that $A_2$  is homogeneous
of degree 2 in $t_3,t_4$. Thus 
\[ F=t_1^2(c_1 t_3^2+c_2 t_3t_4+c_3 t_4^2) +t_1(t_2t_3-t_4^2)(c_4 t_3+c_5 t_4) +(t_2t_3-t_4^2)^2 .\] 
We will show that an irreducible equation $F$ as above, defines a ruled surface.   
Consider a point $(0,1,\mu ^2,\mu)\in D,\ \mu \neq 0,\infty$, then there is a  
$(1,\lambda ,0,0)\in L$ such that the line $\{(s,s\lambda +1,\mu ^2, \mu )|\ s\in \mathbb{P}^1\}$ lies 
on the surface $F=0$. Indeed, substitution in $F$ yields the equation  
\[ s^2(c_1\mu ^4+c_2\mu ^3+c_3\mu ^2)+s^2\lambda \mu ^2(c_4\mu ^2+c_5\mu )+ 
s^2\lambda ^2\mu ^4=0.  \] 
For general constants $c_i$ and general $\mu \neq 0,\infty$, this equations has two solutions for 
$\lambda$. If the equation has for every $\mu$ only one solution for
$\lambda$, then one easily verifies that $F$ is reducible (in fact a square).

Suppose now that $(1,\lambda,0,0)\in L$ is given. The $\mu \neq 0,\infty$ such that the
line $\{(s,s\lambda +1,\mu ^2, \mu )|\ s\in \mathbb{P}^1\}$ lies on $F=0$ are  solutions of the equation
\[\mu ^2(\lambda ^2+c_4\lambda +c_1)+\mu (c_2+\lambda c_5)+c_3=0 .\]
\noindent (a) {\it Suppose $c_3\neq 0$}. If the equation has only one solution for $\mu \neq 0,\infty$, 
then $F$ is easily verified to be a square. The assumption that $F$ is irreducible, implies that 
there are for general $\lambda$ two solutions $\mu$.\\
We conclude that the maps $C\rightarrow D$ and $C\rightarrow L$, given by
$\overline{w}\in C\mapsto \overline{w}\cap H\in D$ and 
$\overline{w}\in C\mapsto\overline{w}\cap L\in L=\overline{e_{12}}$ have both degree 2. 
A further calculation shows that  $P(C)=T_{Gr,\overline{e_{12}}}$,  $\overline{e_{12}}\not \in C$ 
and the vector bundle $B$ has type $-2,-2$. There are still two cases, 
{\it Number} 11 and 12.\\ 

\noindent (b). If $c_3=0$, then $c_2=c_5=0$ is excluded by $F$ is irreducible. Thus there is 
only one solution $\mu$. The maps $C\rightarrow D$ and $C-\rightarrow L=\overline{e_{12}}$ 
have degrees 2 and 1. Further calculation shows
that $P(C)=T_{Gr,\overline{e_{12}}}$, $\overline{e_{12}}\in C$ and the vector bundle
$B$ has type $-1,-3$.  This is {\it Number} 3. \\  

In Rohn's paper only case (a) is considered and this is treated as follows. 
The image $E$ of the morphism $C^{norm} \rightarrow D\times L$ is given by an irreducible  
bi-homogeneous form of bi-degree $(2,2)$. Since $C^{norm}$ has genus 0, the curve $E$ 
has a singular point which is a node or a cusp. The embedding 
$E\subset D\times L\cong \mathbb{P}^1\times \mathbb{P}^1$ can be chosen to be symmetric 
if the field $K$  is algebraically closed. For $K=\mathbb{R}$ one can have a symmetric or an 
anti-symmetric embedding. \\

{\it If $E$ has a node}, then the equation $A$, symmetric in $\lambda ,\mu$, for the embedding is 
written as $a_1\lambda ^2\mu ^2+a_2(\mu ^2\pm \lambda ^2)+2a_3\lambda \mu$, where $\lambda$ 
and $\mu$ are inhomogeneous coordinates for  the rational curves $L$ and  $D$. 
The $\pm$ sign takes care of the real case where one also has   to consider an anti-symmetric 
embedding. The singular point of $E$ corresponds to $\lambda =\mu =0$, which is the point 
$(0,1,0,0)$.   The surface $S$ containing the family of the  lines through the the pairs of points  
$\{(\lambda ,1,0,0), (0,1,\mu ^2,\mu )\}$ satisfying $A(\lambda ,\mu )=0,\ \lambda , \mu \neq 0$ has 
the equation 
\[a_1t_1^2t_3^2+a_2((t_2t_3-t_4^2)^2\pm t_1^2t_4^2)+2a_3t_1t_4(t_2t_3-t_4^2)=0. \]    
 There are various possibilities over $\mathbb{R}$ of the `pinch points' on $L$ and $D$, i.e., the 
ramification points for the two projections $pr_1:E\rightarrow D,\ pr_2:E\rightarrow L$.
\begin{enumerate}  
\item[(i)]  All four are real if $\pm =+$ and $\frac{a_3^2-a_2^2}{a_1a_2}>0$.\\  
Series XIII, no 8, corresponds to this case with additionally  $a_1,a_2>0$.
\item[(ii)] No real ones, if $\pm =+$ and  $\frac{a_3^2-a_2^2}{a_1a_2}<0$.
\item[(iii)] Real on $L$  and not real on $D$ if $\pm =-$ and $a_1>0,\ a_2<0$.
\item[(iv)] Not real on $D$ and real on $L$  if $\pm =-$ and $a_1>0,a_2>0$.
\end{enumerate}

{\it If $E$ has a cusp}, then the equation $A$, symmetric in $\lambda ,\mu$, for the
embedding $E\subset \mathbb{P}^1\times \mathbb{P}^1$ can be normalized (following Rohn) to
\[(\lambda -\mu )^2-2\lambda \mu (\lambda +\mu )+\lambda ^2\mu^2=0.
\mbox{ This leads to the equation}\]
\[t_1^2t_3^2-2t_1t_3(t_1t_4+t_2t_3-t_4^2)+(t_1t_4-t_2t_3+t_4^2)^2=0\mbox{ for } S.\]

\subsubsection{$1^3$}\label{3.2.3}
The line with multiplicity 3 is chosen to be $t_3=t_4=0$. Then the equations have the form 
$t_1A_1+t_2A_2+A_3=0$ with $A_1,A_2,A_3$ homogeneous in $t_3,t_4$ of degree 3,3,4;  
$\gcd(A_1,A_2,A_3)=1$ and $A_1,A_2$ are linearly independent.  Conversely, one  easily verifies 
that the above equation defines a ruled surface of degree 4. 

The pair $(A_1,A_2)$ is unique up to 
taking linear combinations (and linear changes of $t_3,t_4$). In other words
the morphism $f:\mathbb{P}^1\rightarrow \mathbb{P}^1$, given by $(t_3:t_4)\mapsto (A_1:A_2)$, is 
unique and can have degree 3, 2 or 1. In the first case there are many  possibilities for $f$. In the 
second case one can normalize $(A_1,A_2)=d(t_3,t_4)\cdot (t_3^2,t_4^2)$ and in the third case 
$(A_1,A_2)=d(t_3,t_4)\cdot (t_3,t_4)$. The term $A_3$ can be changed into 
$A_3+\ell_1A_1+\ell _2A_2$ with $\ell _1,\ell _2$ homogeneous in $t_3,t_4$ of degrees 1, by 
replacing $t_1,\ t_2$ by $t_1+\ell _1,\ t_2+\ell _2$.
\begin{enumerate}
\item $\gcd(A_1,A_2)=1$ and $A_3=0$.   {\it Number} $1\ a,b,c$. 
\item $\gcd(A_1,A_2)=1$ and $A_3\not \in \{\ell _1A_1+\ell _2A_2\}$. 
{\it Number}  $13\ a,b,c$, XIII 6.
\item $\gcd(A_1,A_2)$ has degree 1. {\it Number} $8,9$.
\item $\gcd(A_1,A_2)$ has degree 2, not a square. $\overline{e_{12}}$ is a node. 
{\it Number} $4$, XIII 7. 
\item $\gcd(A_1,A_2)$ has degree 2 and is a square. $\overline{e_{12}}$ is a cusp.
{\it Number }$5$. 
\end{enumerate} 

\subsubsection{$3^2$}\label{3.2.4}
The twisted cubic curve $TC:=\{(1,\lambda ,\lambda ^2,\lambda ^3)|\ \lambda \in \mathbb{P}^1\}$
is the singular locus of $S$. The homogeneous ideal of $TC$  is generated by the three
 homogeneous forms $X=t_1t_3-t_2^2,\ Y=t_2t_3-t_1t_4,\ Z=t_2t_4-t_3^2$.  
There are two relations $t_3X+t_2Y+t_1Z=t_4X+t_3Y+t_2Z= 0$. The equation $F$ of $S$ is 
homogeneous of degree 4 and lies in the ideal $(X,Y,Z)^2\subset K[t_1,t_2,t_3,t_4]$. 
A computation in the ring $R:=K[\frac{t_2}{t_1},\frac{t_3}{t_1},\frac{t_4}{t_1}] $ shows that the 
element $G:=F(1,\frac{t_2}{t_1},\frac{t_3}{t_1},\frac{t_4}{t_1})$ of total degree $\leq 4$, 
lying in the ideal $(R\frac{X}{t_1^2}+R\frac{Y}{t_1^2})^2\subset R$, is a homogeneous 
polynomial in the terms $\frac{X}{t_1^2},\frac{Y}{t_1^2},\frac{Z}{t_1^2}$ of degree $2$. It follows 
that $F(t_1,t_2,t_3,t_4)=H(X,Y,Z)$, where $H$ is a homogeneous form of degree 2. \\
Consider the morphism $f:\mathbb{P}(V)\setminus TC\rightarrow \mathbb{P}^2$, given by
$(t_1,t_2,t_3,t_4)\mapsto (X,Y,Z)$. The fibres of $f$ are the lines of $\mathbb{P}(V)$ intersecting
$TC$ with multiplicity 2. Thus a fibre is a corde of $TC$ or a tangent line of $TC$.  
Let $H(X,Y,Z)$ be homogeneous of degree 2. Then the closure of the preimage under $f$ of the
subscheme $H=0$ of $\mathbb{P}^2$ is the ruled surface $S$ given by the equation
$F(t_1,t_2,t_3,t_4):=H(X,Y,Z)$. Further $F$ is irreducible and reduced if and only if  
$H=0$ is a conic.  In the sequel we suppose that $\{H=0 \}$ is a conic and we classify the 
possibilities.  The surface with $H=T:=Y^2-4XZ$ is rather special. It consists of 
{\it all tangent lines of $TC$} (see Corollary~\ref{1.2.6}).  For any other conic $H=0$, the intersection 
with $T=0$ has multiplicity 4. {\it In the general case, the intersection of the two conics 
consists of 4 points}. XIII 9, 10.   

Suppose that the intersection has at least one point with multiplicity $>1$. The projective space 
$\mathbb{P}^3$  admits an automorphism which preserves the curve $T=Y^2-4XZ=0$ and brings 
this point to $(X,Y,Z)=(0,0,1)$.  Then $H$ has the form $XZ+aX^2+bXY+cY^2$. One has the 
following cases for the intersection.  
\begin{itemize} 
 \item[(i)] $aX^2+bXY+(c+1/4)Y^2=0$ has two distinct solutions (i.e., $b^2-a(4c+1)\neq 0$) and 
$(c+1/4)\neq 0$. Then the intersection consists  of one point with multiplicity 2 and  two points 
with multiplicity 1. 
\item[(ii)] $aX^2+bXY+(c+1/4)Y^2=0$ has two distinct solutions and $(c+1/4)= 0$. 
Then the intersection consists of one point with multiplicity 3 and one point with multiplicity 1.   
\item[(iii)] $aX^2+bXY+(c+1/4)Y^2=0$ has one solution (i.e. $b^2-a(4c+1)=0$) and $(c+1/4)\neq 0$. 
Then the intersection consists of two points with multiplicity 2. 
\item[(iv)] $aX^2+bXY+(c+1/4)Y^2=0$ has one solution (i.e. $b^2-a(4c+1)=0$), $(c+1/4)=0$ and 
$a\neq 0$. Then the intersection consists of one point with multiplicity 4.  
\end{itemize}   
{\it Thus we found in total six distinct cases for $3^2$}  (compare \cite{Bot}). As we will show below, 
there is a further natural subdivision of these classes. \\   

 Let $C\subset Gr$ be the curve associated to the surface $S_H$ associated to the irreducible 
$H=H(X,Y,Z)$ of   degree 2. The  morphism $C\rightarrow \{H=0\}\subset \mathbb{P}^3$ is clearly 
an isomorphism. Thus $C$ is a nonsingular  rational curve. It is clear that $P(C)$ does not lie 
in two tangent spaces of $Gr$. Moreover, since $C$ is not singular,  one must have 
$\dim P(C)=4$. For the surfaces $S_T$ and $S_H$ with $H=XZ+aX^2-\frac{1}{4}Y^2$ with  
$a\neq 0$ (case (iv)  above), $P(C)$ is not a tangent space.  For the remaining 4 classes   
there are, a priori, now two possibilities:    

\noindent (a) If $P(C)$ is not a tangent space, then $B$ has type $-2,-2$. There
are in total six cases and they fill up  {\it Number } 14.\\  
\noindent (b) If $P(C)$ is a tangent space $T_{Gr,\overline{w}_0}$. Then $B$ has type $-1,-3$. 
{\it Number} $2\ a,b,c$. These cases are explained as follows.\\ 
The line $\overline{w}_0$ coincides with $\overline{e_{12}}$ of the cases $2\ a,b,c$. The image 
$f(\overline{w}_0)$ is the conic given by $H=0$. The possibilities for intersection of $H=0$ with  
$T=0$ reflects the possibilities for the ramification of the degree 3 morphism in $2\ a,b,c$. 
Case $2 \ a$ corresponds to (iii) above;
case  $2 \ b$ to  (i) above; case $2\ c$ to the case where the intersection consists of 4 points.

\subsubsection{ $1^2$} \label{3.2.5}
From the Observations~\ref{1.2.obs} and Subsection~\ref{1.4}, one obtains that the genus of $C$ is 1. 
Further $P(C)$ lies in only one tangent space of $Gr$, say at the point $\overline{e_{12}}$, 
since otherwise the surface $S$ has two skew singular lines.  The morphism 
$C\rightarrow \overline{e_{12}}$, given by 
$\overline{w}\in C\mapsto \overline{w}\cap \overline{e_{12}}$, has degree 2 since $m=2$.
This map has 4 ramification points and we obtain for, say,  $t\neq 0,1,\lambda ,\infty$ two
lines of $C$ through the point $(1,t,0,0)\in \overline{e_{12}}$. The map 
$t\neq 0,1,\lambda ,\infty \mapsto P(t)$, where $P(t)\supset \overline{e_{12}}$ denotes the plane 
through these two lines, has degree 1. We may suppose that $P(t)\cap \{(0,0,*,*)\}=(0,0,1,t)$.   
The equation for $S$ is 
\[t_1^2A+t_2^2B+t_1t_2C+t_1D+t_2E+F;\  A,B,C,D,E,F \mbox{ homogeneous in } t_3,t_4 .\] 
For any point $(a_1,a_2,0,0)\in \overline{e_{12}}$, the plane $a_2t_3-a_1t_4=0$ meets  $S$ in 
$\overline{e_{12}}$ and two lines (or one with multiplicity 2) through $(a_1,a_2,0,0)$. This implies 
that  $t_1^2A(t_3,t_4)+t_2^2B(t_3,t_4)+t_1t_2C(t_3,t_4)$ is a multiple of $(t_4t_1-t_3t_2)^2$ and that  
$t_1D(t_3,t_4)+t_2E(t_3,t_4)$ is divisible by $(t_4t_1-t_3t_2)$. After changing the variables 
$t_1,t_2$   we are reduced to  two possible equations for $S$:    
\[  \ (t_4t_1-t_3t_2)^2+H(t_3,t_4)=0 \mbox{ and }       \ G(t_3,t_4)(t_4t_1-t_3t_2)+H(t_3,t_4)=0 .\]
The line $\overline{e_{12}}$ has multiplicity 3 for the second equation. Thus only the first
equation is possible with $H$ not a square since $S$ is irreducible. Moreover, the ruled surface 
defined by this  equation has $\overline{e_{12}}$ as singular locus if and only if $H$ has no 
multiple factor.  Rohn found an equation of this form, namely        
\[ a(t_3^2\pm t_4^2)+2bt_3^2t_4^2+c(t_4t_2-t_4t_1)^2 =0. \]
The sign $\pm$ distinguishes two classes of real cases.  For $\pm =+$ and $\frac{b}{a}<-1$, 
the four ramification points of $C\rightarrow \overline{e_{12}}$ are real. 
This is {\it Number} 15 and  Series XIII, no 5.\\
     
\noindent  {\it  Remark}. The equation $(t_4t_1-t_3t_2)^2+H(t_3,t_4)=0$ where $H$ has no multiple
 factors,  is valid for any field $K$.  If $K$ is algebraically closed, then $H$ is determined by the 
$j$-invariant of the four zeros of $H$ in $\mathbb{P}^1$. For a general field $K$ there are forms
 for $H$.\\                   

\subsubsection{  $1^2,1^2,\ int $, intersecting lines} \label{3.2.6}  
The two intersecting lines $L_1,L_2$, making up the singular locus of the ruled quartic surface $S$, 
lie in a plane $H$. For $\overline{w}\in C$ and $\overline{w}\neq L_1,L_2$
the intersection $\overline{w}\cap H$ is a point of $L_1\cup L_2$. The induced morphism  
$C^{norm}\rightarrow L_1\cup L_2$ has, say, the line $L_1$ as image. Thus we find a nonconstant
 morphism $f:C^{norm}\rightarrow L_1$ and $P(C)$ lies in the tangent space of $Gr$  
at the point $L_1$. For $q\in L_2$ and $q\not \in L_1$, there is no 
$\overline{w}\in C, \overline{w}\neq L_1,L_2$ with $q\in \overline{w}$. One concludes
that $L_2\in C$.  Moreover $L_2$ is a singular point $s$ of $C$ since $L_2$ belongs to the 
singular locus. In particular, $C$ is a rational curve and $\dim P(C)=3$. If $P(C)$
lies in the tangent space of $Gr$ at another point $M\in Gr$, then one obtains a morphism 
$C\rightarrow M$ by $\overline{w}\mapsto \overline{w}\cap M$. Since $C$ has a singular point, 
this morphism has degree $>1$ and one finds the contradiction that $M$ belongs to the singular 
locus. Thus $P(C)$ lies in a single tangent space.
The rational map $C-\rightarrow L_1$, given by $\overline{w}\mapsto \overline{w}\cap L_1$, 
is well defined at the singular point $s\in C$. Then $f$ has degree $>1$ and its
degree is 2, since $L_1$ has multiplicity 2. Further $L_1\not \in C$, otherwise the multiplicity of 
$L_1$ would be 3. \\
For a suitable basis of $V$ and parametrization of $C^{norm}$, the morphism 
$C^{norm}\rightarrow L_1$ has the form $\overline{w}(t)\mapsto (1,t^2,0,0)$.  Let 
$b:=(b_1,b_2,b_3,b_4)$, with all $b_i\in K[t]$ and $\gcd(b_1,\dots ,b_4)=1$,
be another point of the line $\overline{w}(t)$. By subtracting a multiple of  $(1,t^2,0,0)$ one 
arrives at $\deg b_2\leq 1$.
The Pl\"ucker coordinates of $\overline{w}(t)$  are $(b_2-t^2b_1,b_3,b_4,t^2b_3,t^2b_4,0)$ and 
thus $\deg b_1, \deg b_3,\deg b_4\leq 2$.
The morphism $C\rightarrow \mathbb{P}^1$, by $\overline{w}(t)\mapsto (b_3(t):b_4(t))$, is well 
defined and not constant. Since $C$ is singular, this morphism has degree 2. The corresponding 
degree $2$ morphism $g: C^{norm}\rightarrow \mathbb{P}^1$ factors over $C$. If the singular point 
of $C$ is a cusp for $t=\infty$, then $t=\infty$ is a ramification point and $g$ has the form 
$g(t)=(1:(t+\alpha )^2)$. If the singular point of $C$ is a node, corresponding to $t=\pm 1$, 
then $g(t)=(1:(\frac{at+b}{ct+d})^2)$ also satisfies $g(1)=g(-1)$.
 Hence $g(t)=(1:(\frac{t+\beta}{\beta t+1})^2)$ with $\beta ^2\neq 1$.  

\bigskip
Suppose that $C$ has a cusp, then $(b_3(t),b_4(t))=(1, (t+\alpha )^2)$ and $b_1,b_2$ can be 
normalized to constant multiples of $t$. The condition that $t=\infty$ is a cusp for $C$ implies 
$b_1=0$ and so we arrive at $b=(0,t,1,(t+\alpha )^2)$. The equation reads
\[(t_2t_3-t_1t_4-\alpha ^2t_1t_3+\alpha ^2 t_3^2)^2-t_1t_3t_4(t_3-2\alpha t_1)^2=0.
\mbox{ {\it  Number} } 10\ b .\]

\bigskip
Suppose that $C$ has a node, then $(b_3,b_4)=((\beta t+1)^2,(t+\beta )^2)$ with $\beta ^2\neq 1$. 
For $\beta =0$, one can normalize $b_1,b_2$ to constant multiples of $t$. 
The condition $\overline{w}(1)=\overline{w}(-1)$ implies that $b_1=b_2=ct\neq 0$. The equation reads
\[ c^2t_3t_4(t_3-t_4)^2-(t_1t_4-t_2t_3)^2=0. \mbox{ {\it Number} } 10 \ a .\]
For $\beta ^2\neq 0,1$, one can normalize $b_1,b_2$ to constants and the condition 
$\overline{w}(1)=\overline{w}(-1)$ implies $b_1=b_2=c\neq 0$. The equation reads
\[t_3t_4\{2\beta (t_1-t_2)+c(1-\beta ^2)(t_3-t_4)\}^2 -\{ -t_1(\beta ^2t_3+t_4)+t_2(t_3+\beta ^2t_4\}^2=0 .\]
Again {\it Number } $10\ a$.
Rohn found the two similar equations
\[ at_3t_4^3+(t_1t_4-t_2t_3)^2=0, \mbox{ and }  at_4^4+2bt_3^2t_4^2+(t_1t_4-t_2t_3)^2=0 .\]          

\subsubsection{$1^2,1^2$, skew lines} \label{3.2.7} 
The skew lines can be supposed to be $\overline{e_{12}},\overline{e_{34}}$.
Every monomial in the equation $F$ of $S$ is divisible by one of the terms $t_1^2,t_1t_2,t_2^2$ and 
by one of the terms $t_3^2,t_3t_4,t_4^2$.  Therefore $F$ is bi--homogeneous of degree $(2,2)$ 
and $F$ defines a Zariski closed subset of 
$\overline{e_{12}}\times \overline{e_{34}}\cong \mathbb{P}^1\times \mathbb{P}^1$, which is an 
irreducible curve $E$.
One considers the morphism 
$f:\mathbb{P}(V)\setminus \overline{e_{12}}\cup \overline{e_{34}}\rightarrow 
\mathbb{P}^1\times \mathbb{P}^1$, given  by $(a_1,a_2,a_3,a_4)\mapsto ((a_1,a_2),(a_3,a_4))$. 
Then $S$ is the Zariski closure of $f^{-1}(E)$. The curve $E$ has no singularities since otherwise 
the singular locus of $S$ would contain another line. Thus $E$ is a curve of genus 1. 
One easily sees that $C$ identifies with $E$ and that $P(C)$ lies in the two tangent space of $Gr$ 
at the points $\overline{e_{12}}$ and $\overline{e_{34}}$.   {\it Number} $16$.  

In the above the bases of the two vector spaces $Ke_1+Ke_2$ and $Ke_3+Ke_4$
(or equivalently the parametrization of $\overline{e_{12}}$ and $\overline{e_{34}}$) can be chosen 
in a suitable way. Rohn (see Section 4) shows that for $K=\mathbb{C}$ these bases can be chosen 
such that the equation $F$ becomes symmetric, i.e., $F(t_1,t_2,t_3,t_4)=F(t_3,t_4,t_1,t_2)$. 
For $K=\mathbb{R}$ the results of Rohn are more complicated. These results are essential for the 
understanding of the models in Series XIII, 1,2,3,4 of quartic ruled surfaces with two skew lines of 
singularities.  

\subsubsection{ $1^2,1^2,1^2$}\label{3.2.8} 
Let $L_1,L_2,L_3$ denote the singular lines with multiplicity 2. 

{\it Suppose that the lines $L_1,L_2$ are skew}. 
Then we may suppose $L_1=\overline{e_{12}},\ L_2=\overline{e_{34}}$. From
Lemma~\ref{1.2.10} it follows 
that the equation $F$ of the surface $S$ is bihomogeneous of degree $(2,2)$ in the pairs of 
variables $t_1,t_2$ and $t_3,t_4$. The curve 
$E\subset \overline{e_{12}}\times \overline{e_{34}}\cong\mathbb{P}^1\times \mathbb{P}^1$, defined 
by $F$, has one singular point corresponding to the line $L_3$. This point is a node or a cusp. 
{\it Number} $6,7$. Not in Series XIII. 
The parametrization of $\overline{e_{12}}$ and $\overline{e_{34}}$ can be chosen
(see Section~\ref{section4}) in order to obtain the standard equations of Rohn  
 \[a_1\lambda ^2\mu ^2+a_2(\lambda ^2\pm \mu ^2)+2a_3\lambda \mu =0 \mbox{ and }
\lambda ^2\mu ^2+(\lambda -\mu )^2-2\lambda \mu (\lambda +\mu)=0,\] 
where  $\lambda = \frac{t_2}{t_1},\ \mu =\frac{t_4}{t_3}.$\\

The next case to consider is $L_1\cap L_2,\ L_1\cap L_3, L_2\cap L_3 \neq \emptyset$. The three 
lines cannotlie in a plane $H$ since otherwise the curve $H\cap S$ has degree 6. It follows 
that $L_1\cap L_2\cap L_3$ is one point. We may suppose that $L_1$ is given by 
$t_1=t_2=0$, $L_2$ by $t_1=t_3=0$ and $L_3$ by $t_2=t_3=0$. Every monomial of the 
equation $F$ is divisible by   $t_1^{a_0}t_2^{a_1}$ with $a_0+a_1=2$, by $t_1^{b_0}t_3^{b_1}$ 
with $b_0+b_1=2$ and by $t_2^{c_0}t_3^{c_1}$ with $c_0+c_1=2$. The $t_4$-part of $F$ can only
 be $c\cdot t_1t_2t_3t_4$. If $c=0$, then $F$ defines a cone. Otherwise one can reduce to the 
equation $(t_2^2t_3^2+t_1^2t_3^2+t_1^2t_2^2)+t_4t_1t_1t_3=0$ 
(or equivalently $(t_2t_3+t_1t_3+t_1t_2)^2+t_4t_1t_2t_3=0$).  This equation defines the 
Steiner's Roman surface and the three singular lines are in fact the only lines on this surface.\\
 
\section{ Rohn's symmetric form for bi--degree $(2,2)$}\label{section4}

K.~Rohn proves that over the field $K=\mathbb{C}$,  there is an identification of 
$\overline{e_{12}}\times \overline{e_{34}}$ with $\mathbb{P}^1\times \mathbb{P}^1$ such that the 
equation $F$ of bi--degree $(2,2)$ is symmetric  in the pairs of variables $t_1,t_2$ and 
$t_3,t_4$, i.e.,  $F(t_3,t_4,t_1,t_2)=F(t_1,t_2,t_3,t_4)$. This leads to only a few standard forms 
for $F$.   Over the field $\mathbb{R}$, there are more possibilities. First of all, 
$\overline{e_{12}},\overline{e_{34}}$ can be a pair of conjugated lines over $\mathbb{C}$. 
Secondly, even if $\overline{e_{12}},\ \overline{e_{34}}$ are real lines, then the above 
identification need not be defined over $\mathbb{R}$.    Thirdly, there are various possibilities 
over $\mathbb{R}$ for the ramification points of the two projections $C\rightarrow  \mathbb{P}^1$. 
The models Series XIII, nr. 1,2,3,4 represent some of these cases.
 A `modern version' of this work of Rohn is as follows. \\   

 Consider the closed subset $E$ of  $\mathbb{P}^1\times \mathbb{P}^1$, defined by a 
bi--homogeneous  form $F$ of bi--degree $(2,2)$.  To start we consider the case that $F$ is
 irreducible and $E$ is  nonsingular and  thus $E$ has genus 1. We call the embedding 
$E\subset \mathbb{P}^1\times \mathbb{P}^1$ {\it symmetric} if 
$(p,q)\in E\Rightarrow (q,p)\in E$. 

\begin{theorem}[K. Rohn]\label{4.0.2} 
For a given embedding $E\subset \mathbb{P}^1\times \mathbb{P}^1$ as above,  
there exists an automorphism $f$ of the first factor, such that the new embedding 
$E\subset \mathbb{P}^1\times \mathbb{P}^1
\stackrel{f\times 1}{\rightarrow}\mathbb{P}^1\times \mathbb{P}^1$ is symmetric.
\end{theorem} 

\begin{proof}   The required automorphism $f$ of $\mathbb{P}^1$ has the property 
$(p,q)\in E\Rightarrow  (f^{-1}q,fp)\in E$. In particular,  the morphism 
$C: (p,q)\mapsto (f^{-1}q,fp)$ is an automorphism of $E$ of order 2. We assume that $f$ exists, 
try to find its explicit form and then use this form to produce an $f$ with the required property.  
Some explicit information concerning the automorphisms of order 2 of $E$ is needed. For this 
purpose, we choose a point $e_0\in E$. This makes $E$ into an elliptic curve (and the addition of 
two points $a,b$ is written as $a+b$). Consider the automorphisms  
$\sigma$ and $\tau _a$ (any $a\in E$), given by  $\sigma (p)=-p$ and $\tau _a(p)=p+a$. 
One verifies that the automorphisms of order 2 of $E$ are:\\

\noindent (a) $\sigma \tau _a$ for any point $a$ on $E$,\\
\noindent (b) $\tau_a$ where $a\neq 0$ is a point  of order two on $E$.\\   

Division of $E$ by the action of an element in the first class yields $\mathbb{P}^1$ and division by 
the action of an element in the second class yields an elliptic curve.
Thus the  two projections $pr_i:E \rightarrow \mathbb{P}^1$ correspond to distinct elements 
$\sigma \tau _{a_1}$ and $\sigma \tau _{a_2}$ of order 2 with the property 
$pr_i\circ \sigma \tau _{a_i}=  pr _i$ for $i=1,2$.     

The assumption on $f$ and the definition of 
$C$ are equivalent to  $pr_2(Ce )=f (pr _1(e ))$ for any $e \in E$. Replacing   
$e$ by $\sigma \tau _{a_1}e$ does not change the right hand side.   
Thus $C\sigma \tau _{a_1}e$ is either $Ce$ or  $\sigma \tau _{a_2}Ce$.  The first equality can 
only hold for four elements $e \in E$. Hence the second  equality holds for almost all $e$ and thus 
holds for all $e$. We conclude that   $C\sigma \tau _{a_1}=\sigma \tau _{a_2} C$.    

Suppose that $C=\sigma \tau _c$. The equality  
$\sigma \tau _c\sigma \tau _{a_1}=\sigma \tau _{a_2}\sigma \tau _c$ is equivalent to  
$2c=a_1+a_2$. There are 4 solutions $c$ of this equation.      

Suppose that $C=\tau _c$ with $c$ 
an element of order 2. Then one finds the contradiction  $a_1=a_2$.      

Take $C=\sigma \tau _c$ for some $c$ with  $2c=a_1+a_2$. Define $f$ by the formula   
$f(pr_1(e )):=pr_2(Ce )$. This is well defined because of  
$C\sigma \tau _{a_1}=\sigma \tau _{a_2}C$. It is easily  verified that $f$ is an isomorphism and  
 has the required property.  
 \end{proof}  

Let $E\subset \mathbb{P}^1\times \mathbb{P}^1$ be a symmetric embedding and the homogeneous 
coordinates of the two projective lines are denoted by $x_0,x_1$ and $y_0,y_1$.  Let
$\{p_1,p_2,p_3,p_4\}\subset \mathbb{P}^1$ denote the four ramification points of the projection 
$pr_1:E\rightarrow \mathbb{P}^1$. There is an automorphism $s$ of order two which permutes  
each of the pairs $\{p_1,p_2\}$ and $\{p_3,p_4\}$. The two fixed points of $s$ can be supposed 
to be $0,\infty $ and thus $s$ has the form $s(x_0,x_1)=(x_0,-x_1)$. The four ramification points are
 then $\{(1,\pm d)\}$ and $\{(1,\pm e)\}$. By scaling $(x_0,x_1)\mapsto (x_0,\lambda x_1)$
with $\lambda ^2ed=\pm 1$ we arrive at four ramification points $\{(1,\pm d^{\pm 1}\}$
(with of course $d^4\neq 1$). Write $F=Ay_0^2+By_0y_1+Cy_1^2$. Then the four ramification points
 of $pr_1$ are the zeros of the discriminant $B^2-4AC$ and thus 
$B^2-4AC=x_1^4+bx_1^2x_0^2+x_0^4$ with $b=-(d^2+d^{-2})$.  Then we obtain 
{\it the normal form of K.~Rohn for $F$}: 
\[ a_1(x_0^2y_0^2+x_1^2y_1^2)+a_2(x_0^2y_1^2+x_1^2y_0^2)+2a_3x_0x_1y_0y_1\ \ \ 
\mbox{ with } a_1a_2\neq 0 \] 
\[\mbox{or in Rohn's notation } a_1(\lambda ^2\mu ^2 +1)+a_2(\lambda ^2+\mu ^2)+
2a_3\lambda \mu  \mbox{ with } \lambda =\frac{x_1}{x_0},\ \mu =\frac{y_1}{y_0} ,\]
\[\mbox{and with discriminant }  x_1^4+bx_0^2x_1^2+x_0^4 \mbox{ and } 
b=\frac{a_{1}^2+a_{2}^2-a_{3}^2}{a_{1}a_{2}}\neq \pm 2\ . \ \ \ \ \  \ \        \]     

The above calculations are valid over any algebraically closed field of characteristic $\neq 2$.  
Now we analyze the more complicated situation over the field  $\mathbb{R}$.  Assume that the two 
lines $\overline{e_{12}},\overline{e_{34}}$ and $E$ are defined over $\mathbb{R}$.  Assume 
moreover that $E(\mathbb{ R})$ is not empty (indeed otherwise the real model for the 
corresponding surface has no points).  Fix a real point $e_0$ as the neutral element of $E$. 
The group $E(\mathbb{R})$ is either isomorphic to the circle  $\mathbb{R}/\mathbb{Z}$ 
({\it the  connected case})  or to  $\mathbb{R}/\mathbb{Z}\times \mathbb{Z}/2\mathbb{Z}$ 
({\it the disconnected case}). In the connected case  $E(\mathbb{R})$ has two elements of order 
dividing 2 and in the disconnected case there are 4 such elements. The collection of the real 
automorphisms of order two of $E$ consists of the $ \sigma \tau _a$ with $a\in E(\mathbb{R})$ and 
$\tau _v$ where $v$ is a real point of order 2. 

{\it Now we revisit  the proof of the theorem for the 
case  $K=\mathbb{R}$}.\\

\noindent {\it The connected case}. The fixed points $b$ of $\sigma \tau _{a_1}$ 
(note that$a_1\in E(\mathbb{R})$) are the solutions of $2b =-a_1$. Two of the $b$'s are real.
The other two are complex conjugated. Hence two of the ramification points for
$pr_1:E \rightarrow \mathbb{P}^1$ are real, the other two are complex conjugated. The same
holds for the ramification points of $pr_2:E \rightarrow \mathbb{P}^1$. For the element $c$ with 
$2c=a_1+a_2$ there are two real choices. Thus the real version of the theorem remains valid in 
this case.  Two of the four ramification points are real and the other two are complex conjugated. 
One can normalize such that the ramification points are $\pm d,\ \pm id^{-1}$ and this
leads to Rohn's normal  equation
\[ a_{1}(-x_0^2y_0^2+x_1^2y_1^2)+a_{2}(x_0^2y_1^2+x_1^2y_0^2)+2a_{3}x_0x_1y_0y_1
 \mbox{ with real }a_{1},a_{2},a_{3}. \]

\noindent {\it The disconnected case}. There are 4 real fixed points of 
$\sigma \tau _{a_1}$ if $a_1$ lies in the component of the identity of $E(\mathbb{R})$. 
In the opposite case there are no real solutions of $2b=a_1$. The same holds for 
$\sigma \tau _{a_2}$ and for the solutions of the equation $2c=a_1+a_2$. Hence there are cases 
where no real automorphism$f$ exist. All cases can be listed by:\\
(a) 4 real ramification points for $pr_1$ and for $pr_2$ and 4 real solutions for $c$,\\
(b) no real ramification points for $pr_1$ and  $pr_2$ and 4 real solutions for $c$,\\
(c) 4 real ramification points for $pr_1$, none for $pr_2$ and no real solution for $c$,\\
(d) 4 real ramification points for $pr_2$, none for $pr_1$ and no real solution for $c$.\\

Suppose that $c$ can be chosen to be real. For Rohn's normal form one needs an automorphism 
$s$ permuting each pair $\{p_1,p_2\}$ and
$\{p_3,p_4\}$. One may suppose that each pair is invariant under complex conjugation.
Then the resulting $s$ is also real. For the cases (a) and (b) the standard equation is indeed 
\[ a_{1}(x_0^2y_0^2+x_1^2y_1^2)+a_{2}(x_0^2y_1^2+x_1^2y_0^2)+2a_{3}x_0x_1y_0y_1 \]
\[  \mbox{ and discriminant }  
x_1^4+\frac{a_{1}^2+a_{2}^2-a_{3}^2}{a_{1}a_{2}}x_1^2x_0^2+x_0^4\ ,\]  
with $a_{1},a_{2},a_{3}\in \mathbb{R}$. One easily calculates that 
$\frac{ (a_{1}+a_{2})^2-a_{3}^2}{a_{1}a_{2}}<0$ corresponds to (a) and
$\frac{ (a_{1}+a_{2})^2-a_{3}^2}{a_{1}a_{2}}>0$ corresponds to (b).\\

For the cases (c) and (d) there is no real symmetric normal form for $F$. In case (c) (case (d) is
 similar), Rohn's real normal form could be called {\it half-symmetric}, because of its form 
\[ a_{1}(x_0^2y_0^2-x_1^2y_1^2)-a_{2}(x_0^2y_1^2-x_1^2y_0^2)+2a_{3}x_0x_1y_0y_1 \ . \]
The models 1, 2 and 3 of Series XIII deal with a pair of real skew double lines. In the terminology
of Rohn, a {\it pinch point} is a ramification point for one of the two projections $pr_1,pr_2$ and
 situated on $L_1$ and $L_2$ with the obvious identification of these lines with the two 
$\mathbb{P}^1$'s. Series XIII nr.~1 corresponds to (a), Series XIII nr~2 corresponds to (b) and
 Series XIII no 3 to (c). Rohn also considers the situation where the ruled surface has a 
pair of complex conjugated lines as double lines and produces a standard form and an example, 
namely model 4 of series XIII. 

\bigskip
\noindent  {\it  Rohn's normal form for other curves $E\subset \mathbb{P}^1\times \mathbb{P}^1$ 
of type $(2,2)$}.\\

These normal forms are useful for \S\S~\ref{3.2.2}, \ref{3.2.5}, \ref{3.2.6}, 
and \ref{3.2.8}. There are three cases:\\ 
(a) $E$ is irreducible and  has a node,  \\
(b) $E$ is irreducible and has a cusp and\\
(c) $E$ is reducible or is not reduced.\\

In the following we use the notation and the ideas of the proof of the theorem.\\
\noindent (a). The nonsingular locus of $E$ is, after a choice of a point $e_0$, the group 
$\mathbb{G}_m$. Let $\sigma$ denote the automorphism $x\mapsto -x$ and define $\tau _a(x)=ax$. 
The automorphisms of order 2 of $E$ are $\sigma \tau _a$ (any $a\in \mathbb{G}_m$) and 
$\tau _{-1}$. Dividing $E$ by the action of $\sigma \tau_a$ yields the quotient $\mathbb{P}^1$ and 
dividing by the action of $\tau _{-1}$ yields a rational curve with a double point. Thus the two 
projections $pr_i:E\rightarrow \mathbb{P}^1$ correspond to order two elements $\sigma \tau _{a_i}$ 
for $i=1,2$ with $a_1\neq a_2$. The required  automorphism $C$ of order two
should satisfy $C\sigma \tau _{a_1}=\sigma \tau _{a_2}C$. There are two possibilities for $C$, 
namely $C=\sigma \tau _c$ with $c^2=a_1a_2$. Thus we find a symmetric embedding  
$E\subset \mathbb{P}^1\times \mathbb{P}^1$ for any algebraically closed field of characteristic 
$\neq 2$.\\  

For $\mathbb{R}$ as base field, the situation is more complicated. Suppose that both lines, i.e., the
 two factors $\mathbb{P}^1$, and $E$ are defined over $\mathbb{R}$. We assume that the 
nonsingular locus $E^*$ has a real point  $e_0$. There are two possibilities for  $E^*(\mathbb{R})$, 
namely: (i) $\mathbb{G}_m(\mathbb{R})=\mathbb{R}^*$ and  (ii) $\mathbb{R}/\mathbb{Z}$.\\

In case (i), one has to solve the equation $c^2=a_1a_2$ with $c\in \mathbb{R}^*$. If there is a 
solution, then one has a symmetric embedding $E\rightarrow \mathbb{P}^1\times \mathbb{P}^1$, 
defined over $\mathbb{R}$. In the opposite case, one makes an anti-symmetric embedding 
(by adding some minus signs). The {\it two standard equations} are
\[ a_1\lambda ^2\mu ^2+a_2(\lambda ^2\pm \mu ^2)+2a_3\lambda \mu =0,
\mbox{ with }\lambda =\frac{x_1}{x_0},\ \mu =\frac{y_1}{y_0} .\]

In case (ii), the automorphisms of order two are the maps 
$f_a: x\mapsto -x+a$ (any $a\in \mathbb{R}/\mathbb{Z}$) and $x\mapsto x+1/2$. The last 
automorphism is ruled out because it does not give a $\mathbb{P}^1$ as quotient. Now we have to 
solve $Cf_{a_1}=f_{a_2}C$ for some order two automorphism $C$. The two solutions for $C$ are 
$f_c$ with $2c=a_1+a_2$. There are two solutions for $c\in  \mathbb{R}/\mathbb{Z}$ and therefore 
there is a symmetric embedding. The standard equation is
\[ a_1\lambda ^2\mu ^2+a_2(\lambda ^2 + \mu ^2)+2a_3\lambda \mu =0,\mbox{ with }
\lambda =\frac{x_1}{x_0},\ \mu =\frac{y_1}{y_0} .\]

Finally, there is the possibility that the two lines form a conjugate pair over $\mathbb{R}$.   
[We do not work out the details here.]\\ 
\\
\noindent (b). The nonsingular locus $E^*$ of $E$ is isomorphic to the additive group 
$\mathbb{G}_a$. The automorphism of order two are 
$f_a: x\mapsto -x+a$ (any $a\in \mathbb{G}_a$).  The equation $Cf_{a_1}=f_{a_2}C$ 
(with $a_1\neq a_2$) has a unique solution $C=f_c$ with $2c=a_1+a_2$.
Thus there exists a symmetric embedding $E\subset \mathbb{P}^1\times \mathbb{P}^1$ and this
embedding is unique. The above is valid for any field of characteristic $\neq 2$, because
the group $\mathbb{G}_a$ has no forms.   The standard equation is
\[\lambda ^2\mu ^2+(\lambda -\mu )^2-2\lambda \mu (\lambda +\mu )=0.\]

\noindent (c) For a reducible or nonreduced $E$, Rohn obtains the following standard equations
\[(\lambda +\mu )^2+2a\lambda \mu =0,\ \lambda ^2\mu ^2\pm (\lambda -\mu )^2=0, 
\ (\lambda -\mu )^2=0.\]

\end{document}